\documentclass{etds}

\usepackage{graphicx,amssymb,amsmath}
\graphicspath{{images/}}

\newtheorem{theo}{Theorem}[section]
\newtheorem{prop}{Proposition}[section]
\newtheorem{lemma}{Lemma}[section]
\newtheorem{cor}{Corollary}[section]
 
\newtheorem{definition}{Definition}

\def\N{\mathbb{N}}

\def\Q{\mathbb{Q}}

\def\B{\mathcal{B}}

\def\D{\mathcal{D}}
\def\E{\mathcal{E}}
\def\F{\mathcal{F}}
\def\G{\mathcal{G}}
\def\H{\mathcal{H}}
\def\K{\mathcal{K}}
\def\L{\mathcal{L}}

\def\O{\mathcal{O}}
\def\P{\mathcal{P}}

\def\T{\mathcal{T}}

\def\LEK{\mathcal{K}[\mathcal{L}]}
\def\LLEKK{\mathbb{K}[\mathbb{L}]}
\def\LEKstar{\mathcal{K}^*[\mathcal{L}^*]}
\def\EK{\mathcal{E}\mathcal{K}}

\def\QOP{\mathcal{Q}\mathcal{E}\mathcal{P}} 
\def\QOPstar{\mathcal{Q}\mathcal{E}\mathcal{P}^*} 

\def\OH{\mathcal{O}\mathcal{H}}
\def\CH{\mathcal{C}\mathcal{H}}
\def\COH{\mathcal{C}\mathcal{O}\mathcal{H}}

\def\Nemb{N_\textnormal{emb}}
\def\Nexp{N_\textnormal{exp}}

\def\AA{\mathbb{A}}
\def\BB{\mathbb{B}}
\def\CC{\mathbb{C}}
\def\DD{\mathbb{D}}
\def\EE{\mathbb{E}}
\def\FF{\mathbb{F}}
\def\II{\mathbb{I}}
\def\GG{\mathbb{G}}
\def\HH{\mathbb{H}}
\def\LL{\mathbb{L}}
\def\KK{\mathbb{K}}

\def\NN{\mathbb{N}}
\def\PP{\mathbb{P}}
\def\RR{\mathbb{R}}
\def\SS{\mathbb{S}}
\def\TT{\mathbb{T}}
\def\UU{\mathbb{U}}
\def\VV{\mathbb{V}}
\def\QQ{\mathbb{Q}}
\def\XX{\mathbb{X}}
\def\YY{\mathbb{Y}}

\def\nn{\nonumber}

\def\sse{\subseteq}
\def\restrict{\!\upharpoonright}

\def\into{\hookrightarrow}

\providecommand{\Aut}{\mathop{\rm Aut}\nolimits}
\providecommand{\Age}{\mathop{\rm Age}\nolimits}
\providecommand{\Flim}{\mathop{\rm Flim}\nolimits}
\providecommand{\lev}{\mathop{\rm lev}\nolimits}

\providecommand{\exp}{\mathop{\rm exp}\nolimits}

\newcommand{\expand}[2][*]{({#2}\oplus{#2}^{#1})} 

\newcommand{\bigslant}[2]{{\raisebox{.2em}{$#1$}\left/\raisebox{-.2em}{$#2$}\right.}} 
\newcommand\overeasy{\mathrel{\overset{\makebox[0pt]{\mbox{\normalfont\tiny\sffamily (E)}}}{=}}} 

\newcommand{\cmark}{\checkmark}
\newcommand{\xmark}{\text{\sffamily X}}

\ETDS{1}{2}{1}{2017}

\begin{document}

\runningheads{M. Pawliuk, M. Soki\'c}{Unique ergodicity and Cherlin's classification}

\title{Amenability and unique ergodicity of automorphism groups of countable homogeneous directed graphs}

\author{Micheal Pawliuk\affil{1} and Miodrag Soki\'c}

\address{\affilnum{1}\ University of Toronto)\\
\email{mpawliuk@mail.utoronto.ca}}

\recd{\today}

\begin{abstract}
We study the automorphism groups of countable homogeneous directed graphs (and some additional homogeneous structures) from the point of view of topological dynamics. We determine precisely which of these automorphism groups are amenable (in their natural topologies). For those which are amenable, we determine whether they are uniquely ergodic, leaving unsettled precisely one case (the ``semi-generic" complete multipartite directed graph). We also consider the Hrushovski property. For most of our results we use the various techniques of \cite{AKL12}, suitably generalized to a context in which the universal minimal flow is not necessarily the space of all orders. Negative results concerning amenability rely on constructions of the type considered in \cite{ZUC13}. An additional class of structures (compositions) may be handled directly on the basis of very general principles. The starting point in all cases is the determination of the universal minimal flow for the automorphism group, which in the context of countable homogeneous directed graphs is given in \cite{JLNW14} and the papers cited therein.
\end{abstract}


\section{Introduction}

The Kechris-Pestov-Todorcevic \cite{KPT05} correspondence relates the topological dynamics of automorphism groups of countable homogeneous structures to combinatorial problems related to structural Ramsey theory. Our aim here is to apply this correspondence and the methods of \cite{AKL12,ZUC13}, suitably generalized, to the determination of the amenable automorphism groups associated with countable homogeneous directed graphs, and also to determine which of these are uniquely ergodic; here one specific case remains open.

The KPT correspondence relies on the identification of automorphism groups of countable homogeneous structures with closed subgroups of the full symmetric group on a countable set, with respect to its natural topology, and the more subtle connection due to Fra\"iss\'e between countable homogeneous structures in relational languages and classes of finite structures with certain closure properties, to be described in more detail in \S~\ref{sec:TechPrelim}.

The main points to be taken from \cite{KPT05} relate to the properties of \textit{extreme amenability} and the determination of the \textit{universal minimal flow}. Here a topological group $G$ is said to be extremely amenable if every continuous action on a compact set has a fixed point; a flow is a continuous action on a compact set; and a universal minimal flow is a minimal flow which covers any other minimal flow by a continuous $G$-invariant map.

If one has a countable homogeneous structure $\Gamma$ associated on the one hand with the family $\Age(\Gamma)$ of all finite structures embedding into it, and on the other hand with the group $G = \Aut(\Gamma)$, then according to the KPT correspondence the property of extreme amenability for $G$ is equivalent to the Ramsey property for $\Age(\Gamma)$ (a structural Ramsey theorem). For example, the automorphism group of the randomly ordered random graph is extremely amenable; this is a reformulation of the structural Ramsey theorem for ordered graphs \cite{NR77,NR83,NR89}.

A sharper interpretation of the structural Ramsey theorem for ordered graphs is as a determination of the universal minimal flow for the automorphism group of the (unordered) random graph. Under suitable hypotheses to be reviewed below, the universal minimal flow for an automorphism group can be identified with the set of expansions of the given structure to a category in which the structural Ramsey theorem holds; e.g., to the space of orderings of the random graph in the case at hand.

%
\subsection{Amenability and Unique Ergodicity.}

Using this determination of the universal minimal flow, \cite{AKL12} investigated problems of amenability and unique ergodicity for automorphism groups of countable homogeneous graphs. Here amenability requires, not a fixed point, but a $G$-invariant probability measure, while unique ergodicity requires amenability, but with a \textit{unique} $G$-invariant probability measure. Thus amenability, unique ergodicity, and extreme amenability form a hierarchy of successively stronger properties. In \cite{ZUC13}, Zucker gave examples of countable homogeneous directed graphs (including one tournament) with non-amenable automorphism group.

We will work systematically through the classification of countable homogeneous directed graphs as given by \cite{Cher98}. In \cite{JLNW14} this classification was used to work out the universal minimal flows explicitly. In a number of cases (some already exploited by \cite{ZUC13}) the appropriate Ramsey class is not obtained from expansions by orders, so we will need to reformulate the methods of \cite{AKL12} in a somewhat broader setting.

Imprimitive structures that is, structures carrying nontrivial equivalence relations require some specific attention when forming Ramsey expansions. In the simplest case, where the equivalence relation is a congruence, the analysis may be given in completely general terms. All other cases require individual attention (sometimes in large groups: the main family of examples to be considered is uncountable, but may be treated in a uniform manner).

Evidently the reader will need to know something of the structure of each type of homogeneous directed graph, and the specifics of the expansions, whether by orderings or other additional structure, to classes with the Ramsey property. We will give this in an introductory section which may be used for reference. The specifics of the proof of the Ramsey property are irrelevant here, and even the meaning of the Ramsey property is not very germane, as \cite{AKL12} comes very close to giving a characterization of the properties of amenability and unique ergodicity that we can take as our starting point, once we have reformulated it at the appropriate level of generality.

%
\subsection{Summary of results.}

It will be helpful at the outset to give a chart showing the various results to be obtained. Of course, this chart makes use of notation for specific families of countable homogeneous directed graphs to be discussed in detail a little farther on.

In the classification of homogeneous structures, one usually places the imprimitive (or otherwise degenerate) examples ahead of the primitive ones. We find a different ordering more suitable here. We list a few degenerate examples, then the imprimitive cases involving a congruence (one might say, the highly imprimitive cases), then a few primitive examples that turn out to have non-amenable automorphism groups, and then the more typical cases in which the automorphism groups are amenable and the structures are not particularly degenerate, more or less in order of their structure-theoretic complexity. One of the more exotic structures from our point of view (the semi-generic complete multipartite directed graph) falls somewhere toward the middle of the classification, from that point of view. We also remark that we include in our analysis a property that we have neglected in this introduction, one of several that arise naturally in consideration of the problem of unique ergodicity, and one that is certainly of independent interest.

The following summarizes what is now known, and what remains open, in regards to the amenability of the automorphism groups of structures from Cherlin's classification. For readability, in the table we suppress $\Aut(\Gamma)$ and simply write the structure $\Gamma$. We also suppress the finite structures.

\begin{table}[h]
\centering
\label{table:results}
\begin{tabular}{lllll}
Results \\ \hline
Type                                        & Notation                              & Amenable?                                          & Hrushovski?           & Uniquely ergodic?            														\\ \hline
Composition                                 & $\TT[\II_n], \II_n [\TT]$             & $\Leftrightarrow \TT$                              & $\Leftrightarrow \TT$ & $\Leftrightarrow \TT$, \S~2  														\\
Some weak local orders                      & $\SS(2), \SS(3)$                      & \xmark, \cite{ZUC13}                      	      	& \xmark                & \xmark                       														\\
\begin{tabular}[c]{@{}l@{}}Generic p.o.,\\
				variant, and 2-cover of $\QQ$\end{tabular} & \begin{tabular}[c]{@{}l@{}}$\PP$, \\
																													$\PP(3), \hat{\QQ}$\end{tabular} & \begin{tabular}[c]{@{}l@{}}\xmark, \cite{KSOK12}\\ 
																																																\xmark, \S~5\end{tabular}												& \begin{tabular}[c]{@{}l@{}}\xmark\\ 
																																																																										\xmark\end{tabular} 
																																																																																				& \begin{tabular}[c]{@{}l@{}}\xmark\\ \xmark\end{tabular} \\ \hline
Linear tournament                           & $\QQ$                                 & \cmark, \cite{pes98}						 									 & \xmark                & \cmark, \cite{pes98}                                    \\
Generic tournament                          & $\TT^\omega$                          & \cmark	                                           & $?$                   & \cmark, \cite{AKL12}                                    \\ \hline
2-cover of generic                          & $\hat{\TT^\omega}$                    & \cmark, \S~4.1                                     & $?$                   & \cmark, \S~8.1                                          \\
Generic multipartite                        & $\DD_n$                               & \cmark, \S~4.2                                     & $?$                   & \cmark, \S~8.2                                          \\
Semi-generic multipartite                   & $\SS$                                 & \cmark, \S~4.3                                     & $?$                   & $?$, \S~8.3                                             \\
Generic omitter                             & $\GG_n, \FF(\T)$                      & \cmark, \S~4.4                                     & $?$                   & \cmark, \S~9                            
\end{tabular}
\end{table}

%
\subsection{Amenability.}

To conclude this introduction we will give some indication of the methods used in the analysis of various cases. As the analysis is simpler in the case of amenability, we begin with that.

Whether we deal with amenability or unique ergodicity, our starting point is a prior understanding of the universal minimal flow in concrete combinatorial terms, supplied in our case by \cite{JLNW14}. This can be taken largely as a black box: we begin with a countable homogeneous structure $\Gamma$, and a certain associated countable homogeneous structure $\Gamma^*$ (with the Ramsey property), and rather than considering $\Aut(\Gamma)$ and $\Aut(\Gamma^*)$, one uses the ideas of \cite{KPT05} as developed in \cite{AKL12} to characterize amenability and universal ergodicity in terms of the combinatorics of the two classes $\K := \Age(\Gamma), \K^* := \Age(\Gamma^*)$ consisting respectively of finite structures embedding into $\Gamma$, and finite structures embedding into $\Gamma^*$.

Namely, amenability is equivalent to existence of a ``random $\K^*$-expansion'', which is a function picking out for each $\K$-structure, a probability measure on the set of its expansions to a $\K^*$-structure, with a coherence condition corresponding to embeddings between structures in $\K$. And then, of course, unique ergodicity is equivalent to the uniqueness of this notion of random expansion.

The treatment of amenability is relatively straightforward. To prove amenability, one may first try the uniform measure, according to which all expansions are equally likely. We may mention two of the standard examples.

Example 1.

1.1. If the expansion from $\K$ to $\K^*$ involves adjunction of an additional ordering, then the coherence condition states that the restriction of a random ordering from a given structure to a smaller structure is random, which is probably clear as it stands and at the combinatorial level is a consequence of the fact that there is only one way to order the small structure (ignoring whatever structure is present beyond the order).

1.2. If the structure carries an equivalence relation $E$, then frequently the expansion consists of an ordering in which $E$-classes are convex. This may be thought of as an ordering of the $E$-classes followed by an ordering of each class, and is an example of what we will call a composition.

In this case, if we have two such structures $\AA \leq \BB$, then a random ordering of $\BB$ clearly induces a random ordering of the classes of $\AA$, and then a random order of each class, independently. This proves the amenability of the automorphism group (with the uniform distribution as witness). There are a number of cases in which the expansion from $\Gamma$ to $\Gamma^*$, or from $\K$ to $\K^*$, involves the introduction of an arbitrary order, and the argument of (1.1) proves amenability. There are also a number of cases where we have a congruence $E$ on our structure in which the both amenability and unique ergodicity can be reduced to a treatment of the quotient structure and the structure on each class, in the manner of (1.2).

When we dispose of cases that can be treated by these methods, we are left with several cases which are in fact non-amenable and for which a contradiction is achieved rapidly by examining the meaning of the coherence condition, and some other cases where the appropriate notion of random expansion must be written down, and the corresponding coherence condition checked.

%
\subsection{Unique Ergodicity.}

Now we should say something about the methods used to prove, or disprove, unique ergodicity, which are more sophisticated than those used to check amenability, and frequently involve some computation.

The first known example of a uniquely ergodic group that is neither extremely amenable nor compact was $S_\infty$ as shown in \cite{GW02}. In \cite{AKL12} it was shown that this is a general phenomenon; the automorphism groups of the Fra\"iss\'e limits of (1) the class of $K_n$-free graphs, (2) the class of $r$-uniform hypergraphs, (3) the class of finite metric spaces with distances in a given additive subsemigroup of $\RR^+$, (4) some classes of hypergaphs with forbidden configurations, are all uniquely ergodic. 

The approach in \cite{AKL12} was to examine consistent random expansions of these Fra\"iss\'e classes. First, amenability was established by verifying that the uniform measure was indeed a consistent random expansion. In their cases the reasonable expansions $\mathcal{K}^*$ were usually arbitrary linear orders, which automatically ensures that the uniform measure works, see their Proposition 9.3. Following that, they establish so called quantitative expansion properties, which amount to asymptotic bounds relating to how often (small) ordered structures embed into other (large) ordered structures. These bounds are established by counting the number of order expansions of a structure, and applying the probabilistic McDiarmid inequality. More concrete details will be provided in Section \ref{sec:UE_McD}.

The cases dealt with here are more subtle since the expansions are usually more complicated than arbitrary linear orders. For example, in the case of $\mathbb{D}_n$ the expansions are linear orders that are convex with respect to the equivalence classes of vertices without edges. Quantifiably more complicated is the case of the semi-generic multipartite digraph $\mathbb{S}$, whose precompact expansion is more than just a collection of linear orders. In this case we were unable to establish unique ergodicity, but were able to establish it for slightly tamer expansions.

For the most part, our approach is to use the tools established in \cite{AKL12} suitably generalized to handle precompact expansions. We streamline the methods into a black box theorem, Lemma \ref{lem:QOP_strategy}, which is purely finitary and combinatorial.

\begin{table}[h]
\centering
\label{table:contents}
\begin{tabular}{ll}
Contents                                                        \\ \hline
\S~1 Introduction                                        & \S~7 Unique ergodicity and \\ & McDiarmid's inequality  \\
\S~2 Preliminaries, including composition                & \S~8 The random method                             \\
\begin{tabular}[c]{@{}l@{}} \S~3 Ramsey expansions of \\ \hspace{0.5cm}countable homogeneous \end{tabular} & \S~9 The hypergraph method						  						\\
\S~4 Amenability                                         & \S~10 Conclusions and open questions \\
\S~5 Failures of amenability                             & \S~11 Appendix \\
\S~6 The Hrushovski property                       			                             
\end{tabular}
\end{table}


\section{Preliminaries, including composition}
\label{sec:TechPrelim}

Now we will describe the mathematical objects and notions that we will be using. Sections 2.1-2.6 are intended to be used as a reference, and Sections 2.7-2.10 are a discussion of compositions. 

%
\subsection{Amenability.}

Let $G$ be a topological group. A continuous action of $G$ on a compact Hausdorff space is called a \textbf{$G$-flow}. A $G$-flow is \textbf{minimal} if the orbit of every point is dense. If every $G$-flow has a $G$-invariant Borel probability measure, then we say that $G$ is \textbf{amenable}. We say that $G$ is \textbf{uniquely ergodic} if every minimal $G$-flow has a unique $G$-invariant Borel probability measure. We will go into more depth about various equivalent versions of amenability in Chapter 2.

Throughout, we consider amenability and unique ergodicity for a collection of automorphism groups of countable structures related to directed graphs. These groups are not locally compact and they are not discrete, but they are non-Archimedian Polish groups, see \cite{BK96} for more details.

%
\subsection{Fra\"iss\'e Classes and Structures.}

Let $\AA$ and $\BB$ be given structures. If there is an embedding from $\AA$ into $\BB$ then we write $\AA\into \BB$, if $\AA$ is a substructure of $\BB$ then we write $\AA \leq \BB$, and if $\AA$ and $\BB$ are isomorphic, then we write $\AA \cong \BB$. We write $\binom{\BB}{\AA} = \{\CC \leq \BB : \CC \cong \AA\}$. We say that a given structure is \textbf{locally finite} if each of its finitely generated substructures are finite. We denote by $\Age(\AA)$ the collection of all finite substructures of $\AA$. A structure $\AA$ is \textbf{ultrahomogeneous} if every isomorphism between two finite substructures can be extended to an automorphism of $\AA$. We say that $\AA$ is a \textbf{Fra\"iss\'e structure} if it is countably infinite, locally finite and ultrahomogeneous.

Let $L$ be a signature and let $\K$ be a class of finite structures in $L$. Then $\K$ satisfies the:

\begin{table}[h]
\begin{tabular}{lp{0.8\textwidth}}
(\textbf{HP}) 	& \textbf{Hereditary Property}, if whenever $\AA \into \BB$ and $\BB \in \K$, then $\AA \in \K$. \\
(\textbf{JEP}) 	& \textbf{Joint Embedding Property}, if for all $ \AA, \BB \in \K$ there is a $\CC \in \K$ such that $\AA~\into~\CC$ and $\BB~\into~\CC$. \\
(\textbf{AP}) 	& \textbf{Amalgamation Property}, if for all $\AA,\BB,\CC \in \K$ and all embeddings $f: \AA \rightarrow \BB$ and $g: \AA \rightarrow \CC$ there is a $\DD \in \KK$ and embeddings $\overline{f}: \BB \rightarrow \DD$ and $\overline{g}: \CC \rightarrow \DD$ with $\overline{f} \circ f = \overline{g} \circ g$. \\
(\textbf{SAP}) 	& \textbf{Strong Amalgamation Property}, if in addition to \textbf{AP} we have $\overline{f}(\BB) \cap \overline{g}(\CC) = \overline{f} \circ f (\AA)$.
\end{tabular}
\end{table}

We say that $\K$ is a \textbf{Fra\"iss\'e class} if it satisfies \textbf{HP}, \textbf{JEP}, \textbf{AP}, contains finite structures of arbitrarily large finite cardinality, and only countably many different isomorphism types. If $\KK$ is a Fra\"iss\'e structure then its $\Age(\KK)$ is a Fra\"iss\'e class. Given a Fra\"iss\'e class $\K$ we have its Fra\"iss\'e limit $\text{Flim}(\K)$, which is a Fra\"iss\'e structure and is unique up to isomorphism. In this way there is a 1-1 correspondence between Fra\"iss\'e classes and Fra\"iss\'e structures. For more details, see \cite{HO93}.

We consider a structure as a tuple $\AA = (A, \{R_i^A\}_{i \in I}, \{f_j^A\}_{j \in J})$ where $A$ is the underlying set of the structure, $R_i^A$ is the interpretation of the relational symbol in $\AA$ and $f_j^A$ is the interpretation of the functional symbol in $\AA$ for all $i \in I$ and all $j \in J$. If $J = \emptyset$ then we say that the structure is \textbf{relational} (or that the signature is relational). All of the structures studied within are relational, so for ease of notation we will appropriate $J$ to also serve as an index set for a collection of relations.

In particular, we consider a directed graph (digraph) as a structure in the binary relational signature $\{\rightarrow\}$. The symbol $\rightarrow$ is always interpreted as an irreflexive and asymmetric relation. For a directed graph $(A, \rightarrow^A)$ we sometimes use the symbol $\perp^A$ to denote $\neg(x \rightarrow^A y \vee y \rightarrow^A x)$. A \textbf{tournament} is a digraph $(A, \rightarrow^A)$ with the property that for every $x \neq y \in A$ we have either $x \rightarrow^A y$ or $y \rightarrow^A x$ (but not both).

In general, we will use the following typefaces: $\L, \K$ for classes, $\AA, \BB, \CC, \KK$ for structures, $A,B,C$ for universes (or underlying sets) of structures (with $L$ being reserved for the signature of a class) and $a,b,c$ for elements of underlying sets. Occasionally we will use $a,b$ for natural numbers that index the number of equivalence classes in $\AA$ and $\BB$.

%
\subsection{Reducts and the expansion property.}

Let $L \sse L^*$ be given signatures. Let $\K$ be a class of structures in $L$ and let $\K^*$ be a class of structures in $L^*$. If $\AA^* \in \K^*$ then we denote by $\AA^* \vert L$ the structure in $\K$ obtained by dropping the interpretations of the symbols in $L^* \setminus L$ in $\AA^*$, and define $\K^* \vert L := \{\AA^* \vert L : \AA^* \in \K^*\}$. We say that $\K^*$ is a \textbf{precompact expansion} of $\K$ provided that $\forall \AA \in \K$ there are only finitely many $\AA^* \in \K^*$ such that $\AA^* \vert L = \AA$.

We say that $\K^*$ satisfies \textbf{the expansion property (EP)} (with respect to $\K$) if $\K^* \vert L = \K$ and for every $\AA \in \K$ there is a $\BB \in \K$ such that for every $\AA^*, \BB^* \in \K^*$ with $\AA^* \vert L = \AA$ and $\BB^* \vert L = \BB$ we have $\AA^* \into \BB^*$.

We say that $\K^*$ is a \textbf{reasonable expansion} of $\K$ provided that it is a precompact expansion and $\forall \AA, \BB \in \K$, for every embedding $\pi : \AA \longrightarrow \BB$, $\forall \AA^* \in \K^*$ with $\AA^* \vert L = \AA$, there is a $\BB^* \in \K^*$ such that $\BB^* \vert L = \BB$ and $\pi$ is also an embedding of $\AA^*$ into $\BB^*$.

%
\subsection{Ramsey property.}

We say that the class $\K$ satisfies the \textbf{Ramsey Property (RP)} (or is a \textbf{Ramsey class}) if for every (small) $\AA \in \K$ and every (medium) $\BB \in \K$ and every $r \in \N$ there is a (large) $\CC \in \K$ such that for every colouring $c : \binom{\CC}{\AA} \rightarrow \{1, \ldots, r\}$ there is a $\overline{\BB} \in \binom{\CC}{\BB}$ such that $c \restrict \binom{\overline{\BB}}{\AA}$ is a constant. We denote this using the arrow notation: 
\[
	\CC \longrightarrow (\BB)_r^\AA.
\]

Another notion that will always appear with \textbf{RP} is rigidity. We say that a given structure is \textbf{rigid} if it has no nontrivial automorphisms.

%
\subsection{Consistent random expansions.}
\label{sec:cons_rand_exp}

Let $L \sse L^*$ be given signatures. Let $\K$ and $\K^*$ be classes of structures in $L$ and $L^*$ respectively such that $\K^*\vert L = \K$. If $\BB \in \K^*$ then we write $\BB^* = \expand{\AA}$ where $\AA = \BB^* \vert L$ and $\AA^* = \BB^* \vert (L^* \setminus L)$. Colloquially, ``$\AA$ is the old stuff, and $\AA^*$ is the new stuff'' when we use the representation $\expand{\AA}$. For $\AA \in \K$ we denote by $\mu_\AA$ a measure on the set 
\[
	\K^*(\AA) := \{\AA^* : \expand{\AA} \in \K^*\}.
\]
We will also have need for the related quantity
\[
	\#(\AA) := \vert \K^*(\AA) \vert,
\]
which is the number of expansions of $\AA$ in $\K^*$. Let $\AA \leq \BB$ be structures in $\K$ and let $\expand{\AA} \in \K^*$. Then we write:
\[
	\#_{\K^*} (\AA^*, \BB) := \left\vert \{ \BB^* : \expand{\AA} \leq \expand{\BB} \in \K^* \}\right\vert,
\]
which is the number of expansions of $\BB$ in $\K^*$ that extend $\expand{\AA}$. If there is no confusion then we write $\#(\AA^*,\BB)$, or occasionally we will write $\#_{\AA,\BB}(\AA^*)$.

We say that the collection $\{\mu_\AA : \AA \in \K\}$ is a \textbf{consistent random $\K^*$-expansion on $\K$} (when it is clear from context we suppress the reference to $\K^*$) if we have:
\begin{itemize}
\item[(\textbf{P}):] Each $\mu_\AA$ is a probability measure on $\K^*(\AA)$.
\item[(\textbf{E}):] Whenever $\varphi: \AA \longrightarrow \BB$ is an embedding, and $\expand{\AA} \in \K^*$, we have $\displaystyle{\mu_\AA (\{\AA^*\}) = \sum \{\mu_\BB (\{\BB^*\}) : \varphi \text{ embeds } \expand{\AA} \text{ into } \expand{\BB}\}.}$
\end{itemize}

\noindent We have \textbf{P}robability measures and \textbf{E}xtension properties. When it is clear from context we shall refer to a consistent random expansion as $(\mu_\AA)$, with no reference to $\K$. In the special case that $\varphi$ in (\textbf{E}) is an isomorphism we get that $\varphi_* \mu_\AA = \mu_\BB$, where $\varphi_* \mu_\AA$ is the push forward measure; call this (\textbf{I}) for \textbf{I}somorphism invariance. We will reference it explicitly later on.

We assume that for $\AA \cong \BB$ in $\K$ and $\expand{\AA} \in \K^*$ we have $\BB^*$ such that $\expand{\AA} \cong \expand{\BB}$.

Let $L^* \setminus L$ be a relational signature, and let $\K$ and $\K^*$ be Fra\"iss\'e classes such that $\K^*$ is a reasonable expansion of $\K$. Then we say that $(\K, \K^*)$ is an \textbf{excellent pair} if:
\begin{enumerate}
	\item $\K^*$ is a Ramsey class of rigid structures, and
	\item $\K^*$ satisfies the expansion property relative to $\K$.
\end{enumerate}

%
\subsection{Amenability via expansions.}

The following is the key equivalence used to show amenability and non-amenability of the automorphism group of a Fra\"iss\'e structure. The version that appears as Proposition 9.2 in \cite{AKL12} is a special case of what we state, and the proof of this version is analogous. The arguments in the proof are standard and the proof uses the Carath\'eodory Extension Theorem, so the proof is omitted.

\begin{prop}[The Key Equivalence]
\label{threeone} 

Let $(\K, \K^*)$ be an excellent pair. Then:
	\begin{enumerate}
		\item $\Aut(\Flim(\K))$ is amenable iff $\K$ has a consistent random $\K^*$-expansion.
		\item $\Aut(\Flim(\K))$ is uniquely ergodic iff $\K$ has a unique consistent random $\K^*$-expansion.
	\end{enumerate} 
\end{prop}

We remark that a consistent random expansion $(\mu_\AA)$ cannot be degenerate, which means that when $\expand{\AA} \in \K^*$, we have $\mu_\AA(\{\AA^*\}) \neq 0$. Otherwise, since $\Flim(\K)$ is separable, a degenerate measure for $\AA^*$ would give us a countable cover of the universal minimal flow by open sets each with measure 0. For more details, see the proof of \cite[Proposition 2.1]{KSOK12}.

%
\subsection{Compositions - $\EK$ and $\LEK$.}
\label{sec:EK_LEK_def}

Now we introduce the composition class $\LEK$.

We mainly focus on the quotient structure $\LEK$, but to introduce it we first mention the class $\EK$, which is used to define $\LEK$. Intuitively, a structure in $\LEK$ is a (horizontal) structure $\KK \in \K$ and associated to each point $k \in \KK$ is a (vertical, possibly different) $\LL_k \in \L$. The $\K$-relations of elements in different $\LL_k$ ``columns'' are given by looking at the $\K$-relations of the corresponding $k \in \KK$. In this way, if you ``quotient out'' by the equivalence relation of being in the same $\LL_k$ column, then you get $\KK$.

Alternatively, one can think of a structure in $\LEK$ as taking a structure $\KK \in \K$ then ``blowing-up'' each of its points $k \in \KK$ to a structure $\LL_k \in \L$.

\begin{definition}[$\EK$] Let $\K$ be a class of structures in $L_I := \{R_i : i \in I\}$, a relational signature where each $R_i$ has arity $n_i$, and let $\sim$ be a binary relational symbol such that $\sim \, \notin L_I$.

We denote by $\EK$ the class of relational structures of the form 
\[
	\AA = (A, \{R^A_i\}_{i \in I}, \sim^A)
\]
with the properties:
	\begin{enumerate}
		\item $\sim^A$ is an equivalence relation on $A$ with equivalence classes denoted by $[a]_{\sim^A}$.
		\item For $i \in I$ and $x_1, \ldots, x_{n_i}, y_1, \ldots, y_{n_i} \in A$ with $x_j \sim^A y_j$ (for all $j \leq n_i$) we have 
			\[
				R^A_i (x_1, \ldots, x_{n_i}) \Leftrightarrow R^A_i (y_1, \ldots, y_{n_i}).
			\]
		\item Let $\bigslant{A}{\sim^A} := \{[a]_{\sim^A} : a \in A\}$ be the set of equivalence classes. Let $R_i^\bigslant{A}{\sim^A}$, for $i \in I$, be the relation defined on the set $\bigslant{A}{\sim^A}$ according to (2) with 
			\[
				R_i^\bigslant{A}{\sim^A}([a_1]_{\sim^A}, \ldots, [a_{n_i}]_{\sim^A}) \Leftrightarrow R_i^A(x_1, \ldots, x_{n_i})
			\]
		 where $a_j \sim^A x_j$ for all $j \leq n_i$. Then we have 
			\[
				\bigslant{\AA}{\sim^A} := (\bigslant{A}{\sim^A}, \{R_i^\bigslant{A}{\sim^A}\}_{i \in I}) \in \K.
			\]
	\end{enumerate}
\end{definition}

\begin{definition}[$\LEK$] Let $L_I := \{R_i : i \in I\}$ and $L_J := \{R_j : j \in J\}$ be disjoint relational signatures and let $\sim$ be a binary relational symbol such that $\sim \, \notin L_I \cup L_J$. Let $\K$ and $\L$ be classes of relational structures in $L_I$ and $L_J$ respectively. 

We denote by $\LEK$ the class of relational structures of the form 
\[
	\AA = (A, \{R^A_i\}_{i \in I}, \{R^A_j\}_{j \in J}, \sim^A)
\]
with the properties:
	\begin{enumerate}
		\item $\AA \vert (L_I \cup \{\sim\}) \in \EK$.
		\item For $j \in J$ and $x_1, \ldots, x_{n_j} \in A$ we have 
			\[
				R^A_j (x_1, \ldots, x_{n_j}) \Rightarrow [x_1]_{\sim^A} = \ldots = [x_{n_j}]_{\sim^A}.
			\]
		\item For $a \in A$ we have 
			\[	
				([a]_{\sim^A}, \{R_j^A \cap ([a]_{\sim^A})^{n_j}\}_{j \in J}) \in \L.
			\]
	\end{enumerate}
\end{definition}

%
\subsection{Expansions of $\LEK$.}
\label{sec:LEK_exp}

Let $L_I^* \supset L_I$ and $L_J^* \supset L_J$ be relational signatures such that $L_I^* \cap L_J^* = \emptyset$ and $\sim \, \notin L_I^* \cup L_J^*$. If $\K^*$ and $\L^*$ are expansions of the classes $\K$ and $\L$ such that $\K^* \vert L_I = \K$ and $\L^* \vert L_J = \L$ then we have that 
\[
	\left(\LEKstar\right) \vert ( L_I \cup L_J \cup \{\sim\} ) = \LEK.
\]

Let $\AA \in \LEK$ be the finite structure which has $A_1, A_2, \ldots, A_a$ as its $\sim^A$-equivalence classes. Let $\AA_1, \AA_2, \ldots, \AA_a$ be structures in $\L$ which are placed on $A_1, \ldots, A_a$ respectively and let $\BB \in \K$ be the structure given by representatives of the equivalence classes. Then we write $\AA = (\BB : \AA_1, \ldots, \AA_a)$.

Similarly, an expansion $\expand{\AA} \in \LEKstar$ of $\AA \in \LEK$ is given by the structures $\expand{\AA_1}, \ldots, \expand{\AA_a} \in \L^*$ and $\expand{\BB} \in \K^*$. So we write
\[
 \expand{\AA} = (\expand{\BB} : \expand{\AA_1}, \ldots, \expand{\AA_a})
\]
or, if there is no confusion
\[
 \AA^* = (\BB^* : \AA_1^*, \ldots, \AA_a^*).
\]

%
\subsection{Excellent pair proposition.}
\label{sec:LEK_excellent}

The following technical proposition ensures that $(\LEK, \LEKstar)$ is an excellent pair, thus we may apply Proposition \ref{threeone} to verify amenability of $\Aut(\Flim(\LEK))$.

\begin{prop}\label{ExPair_LEK} Let $L_I^* \supset L_I$ and $L_J^* \supset L_J$ be relational signatures such that $L_I^* \cap L_J^* = \emptyset$ and let $\sim$ be a binary relational symbol such that $\sim \, \notin L_I^* \cup L_J^*$. Let $\K,\K^*, \L$ and $\L^*$ be classes of finite relational structures in $L_I, L_I^*, L_J$ and $L_J^*$ respectively. Let $\K^* \vert L_I = \K$ and $\L^* \vert L_J = \L$. Then we have:

	\begin{enumerate}
		\item If $\L$ and $\K$ are Ramsey classes of rigid structures then $\LEK$ is a Ramsey class of rigid structures.
		\item If $\L^*$ satisfies \textbf{EP} with respect to $\L$ and $\K^*$ satisfies \textbf{EP} with respect to $\K$, then $\LEKstar$ satisfies \textbf{EP} with respect to $\LEK$. 
	\end{enumerate}
\end{prop}

\proc{Proof.} This follows by simple modifications of the proofs of Theorem 4.4 and Proposition 5.2 in \cite{SOK13}.
\ep
\medbreak

%
\subsection{$\LEK$.}
\label{sec:LEK}

The following theorem is the main result of this section. Establishing this theorem was the genesis of this larger project, and after it was established we expanded our aims to the other digraphs on Cherlin's classification. The forward implication in each of the parts is a straightforward, if tedious combinatorial verification. It can also be derived at the level of topological groups. The converse is more subtle so we include its proof.

\begin{theo}\label{A_LEK} Let $(\L,\L^*)$ and $(\K,\K^*)$ be excellent pairs of classes of finite structures in distinct signatures. Then we have:
	\begin{enumerate}
		\item $\Aut(\Flim(\LEK))$ is amenable iff $\Aut(\Flim(\L))$ and $\Aut(\Flim(\K))$ are amenable. \label{A_LEK_one}
		\item $\Aut(\Flim(\LEK))$ is uniquely ergodic iff $\Aut(\Flim(\L))$ and $\Aut(\Flim(\K))$ are uniquely ergodic.
	\end{enumerate}
\end{theo}

As an (almost) immediate corollary we get the unique ergodicity of $\Aut(\TT[\II_n])$ and $\Aut(\II_n[\TT])$, which are both part of Cherlin's classification.

The proof of this is broken up into five not entirely independent parts. The consistent random expansions presented in the amenability proofs will be used in the unique ergodicity proofs. Moreover, in order to not overly repeat ourselves, detailed proofs that some maps are actually consistent random expansions will only appear in the amenability proofs. These proofs are ``direct'' in the sense that they do not rely on heavy machinery. The essential claim in these proofs is that an expansion in the composition class is a composition of expansions, but there are details that need to be checked.

In what follows the summations will always range over $\K^*$-expansions of fixed structures in $\K$. 
\proc{Proof of (i), $\Leftarrow$.}

Assume that $\Aut(\Flim(\L))$ and $\Aut(\Flim(\K))$ are amenable. Then by Proposition \ref{threeone} there are consistent random expansions $\nu$ and $\mu$ on $\L$ and $\K$ respectively. We will define a consistent random expansion $\nu \otimes \mu$ on $\LEK$.

Let $\SS = (\AA : \SS_1, \ldots, \SS_a) \in \LEK$ and $\expand{\SS} \in \LEKstar$ such that $\SS^* = (\AA^* : \SS_1^*, \ldots, \SS_a^*)$.

Define
	\[
		(\nu \otimes \mu)_\SS (\{\SS^*\}) := \mu_\AA (\{\AA^*\}) \cdot \prod_{i=1}^a \nu_{\SS_i} (\{\SS_i^*\}),
	\]
and we check the conditions for being a consistent random expansion.

$(\textbf{P})$ Observe that
\begin{align*}
	&\space \sum \left\{(\nu \otimes \mu)_\SS (\{\SS^*\}) : \expand{\SS} \in \LEKstar\right\} 																															\\
	&= \sum \left\{\mu_\AA (\{\AA^*\}) \cdot \prod_{i=1}^a \nu_{\SS_i} (\{\SS_i^*\}) : \expand{\AA}\in\K^*, \expand{\SS_i} \in \L^*, \forall i \leq a\right\}	\\
	&= \left(\sum_{\expand{\AA}\in\K^*} \mu_\AA (\{\AA^*\}) \right) \cdot 
							\sum \left\{\prod_{i=1}^a \nu_{\SS_i} (\{\SS_i^*\}) : \expand{\SS_i} \in \L^*, \forall i \leq a\right\} 																	\\
	&= 1 																				\cdot \prod_{i=1}^a \left\{\sum \nu_{\SS_i} (\{\SS_i^*\}) : \expand{\SS_i} \in \L^* \right\}					\\
	&= 1 \cdot \prod_{i=1}^a 1 = 1,																																																			
\end{align*}
where the third equality follows from (\textbf{P}) on $\mu$, and the fourth equality follows from (\textbf{P}) on $\nu$. The second equality follows from a basic fact about sums and products.

(\textbf{E}) Let $\SS = (\AA : \SS_1, \ldots, \SS_a)$ and $\TT = (\BB : \TT_1, \ldots, \TT_b)$ be structures in $\LEK$ such that $\SS \leq \TT$. Let $\SS^* = (\AA^* : \SS_1^*, \ldots, \SS_a^*)$ be such that $\expand{\SS}\in \LEKstar$. Since $\SS \leq \TT$ there is an $I \sse \{1, \ldots, b\}$ such that $\SS_i \leq \TT_i$ if and only if $i \in I$. We take $J = \{1, \ldots, b\}\setminus I$. Then
\begin{align*}
	&  \sum \{(\nu \otimes \mu)_\TT (\{\TT^*\}) : \expand{\SS} \leq \expand{\TT} \in \LEKstar \}	\\
	&= \sum \left\{ \mu_\BB(\{\BB^*\}) \cdot \prod_{i=1}^b \nu_{\TT_i}(\{\TT_i^*\}) :
					\begin{aligned}
						&\scriptstyle{\expand{\AA} \leq \expand{\BB} \in\K^* }\\
	 					&\scriptstyle{\expand{\SS_i} \leq \expand{\TT_i}\in\L^*, \text{ for } i \in I} \\
	 					&\scriptstyle{\expand{\TT_i} \in \L^*, \text{ for } i \in J}
					\end{aligned}
				\right\} \\
	&=  \left( \sum \{\mu_\BB (\{\BB^*\}) : \expand{\AA} \leq \expand{\BB} \in \K^*\} \right)\\
	&\hspace{0.7cm}\cdot \sum \left\{ \prod_{i=1}^b \nu_{\TT_i} (\{\TT_i^*\}) : 
								\begin{aligned}	
									&\scriptstyle{\expand{\SS_i} \leq \expand{\TT_i} \in \L^*, \text{ for } i \in I} \\
									&\scriptstyle{\expand{\TT_i}\in\L^*, \text{ for } i \in J}
								\end{aligned} 
							\right\} 	\\
	&= \mu_\AA (\{\AA^*\}) 
			\cdot \prod_{i\in I} \left\{ \sum \nu_{\TT_i}(\{\TT_i^*\}) : \expand{\SS_i} \leq \expand{\TT_i} \in \L^* \right\} \\
			&\hspace{0.7cm}\cdot \prod_{j\in J} \left\{ \sum \nu_{\TT_j}(\{\TT_j^*\}) : \expand{\TT_j} \in \L^* \right\} \\
	&= \mu_\AA (\{\AA^*\}) \cdot \prod_{i \in I} \nu_{\SS_i}(\{\SS_i^*\}) \cdot \prod_{j \in J} 1
		= \mu_\AA (\{\AA^*\}) \cdot \prod_{i =1}^a \nu_{\SS_i}(\{\SS_i^*\})
		= (\nu \otimes \mu)_\SS (\{\SS^*\}),
\end{align*}
where the third equality comes from (\textbf{E}) of $\mu$, the fourth equality comes from (\textbf{E}) of $\nu$ and the second equality uses the same basic fact about sums and products.

This completes the verification that $\nu \otimes \mu$ is a consistent random expansion, and by Proposition \ref{threeone} we have that $\Aut(\Flim(\LEK))$ is amenable.
\ep\medbreak

\proc{Proof of (i), $\Rightarrow \Aut(\Flim(\K))$ is amenable.} 

Assume that $\Aut(\Flim(\LEK))$ is amenable. By Proposition (\ref{threeone}) there is a consistent random expansion $\rho$ on $\LEK$. We will show that $\Aut(\Flim(\K))$ is amenable.

Let $\KK \in \K$ and $\expand{\KK}\in\K^*$. Let $\SS = (\KK : \SS_1, \ldots, \SS_a) \in \LEK$. Consider
\[
	\mu_{\KK,\SS}(\{\KK^*\}) := \sum \{\rho_\SS(\{\SS^*\}) : \SS^* = (\KK^* : \SS_1^*, \ldots, \SS_a^*), \expand{\SS}\in \LEKstar\}.
\]

First, we show that $\mu_{\KK,\SS}$ is independent of our choice of $\SS$. Let $\TT = (\KK : \TT_1, \ldots, \TT_a)$ also be a structure in $\LEK$. If $\SS \leq \TT$, then
\begin{align*}
	\mu_{\KK, \TT}(\{\KK^*\}) &= \sum \{\rho_\TT (\{\TT^*\}) : \TT^* = (\TT_1^*, \ldots, \TT_a^*: \KK), \expand{\TT}\in\LEKstar\} 	\\
	&= \sum_{\substack{\SS^* = (\KK^* : \SS_1^*, \ldots, \SS_a^*) \\ \expand{\SS}\in \LEKstar }} 
			\sum\{\rho_\TT (\{\TT^*\}) : \expand{\SS} \leq \expand{\TT} \in \LEKstar\} 																					\\
	&= \sum_{\substack{\SS^* = (\KK^* : \SS_1^*, \ldots, \SS_a^*) \\ \expand{\SS}\in \LEKstar }} \rho_\SS (\{\SS^*\}) 			\\
	&= \mu_{\KK,\SS}(\{\KK^*\})
\end{align*}

Where the third equality follows from (\textbf{E}) of $\rho$. If $\SS$ is not a substructure of $\TT$, then by \textbf{JEP} for $\LEK$ there is an $\RR \in \LEK$ such that $\SS \leq \RR$ and $\TT \leq \RR$. So by the above we have:
\[
	\mu_{\KK, \SS}(\{\KK^*\}) = \mu_{\KK, \RR}(\{\KK^*\}) = \mu_{\KK, \TT}(\{\KK^*\})
\]

Therefore $\mu_{\KK, \SS}(\{\KK^*\})$ is independent of the choice of structure $\SS$, so without ambiguity, we write
\[
	\mu_\KK (\{\KK^*\}) := \mu_{\KK, \SS}(\KK^*)
\]
where $\SS = (\KK : \SS_1, \ldots, \SS_a) \in \LEK$.

Now we check that $\mu_\KK$ is a consistent random expansion for $\K$.

(\textbf{P}) Fix any $\SS = (\SS_1, \ldots, \SS_a: \KK) \in \LEK$. Observe that
\begin{align*}
	&  \sum \{\mu_\KK (\{\KK^*\}) : \expand{\KK} \in \K^*\} 																																						\\
	&= \sum \{\mu_{\KK,\SS}(\{\KK^*\}) : \expand{\KK} \in \K^*\} 																																			\\
	&= \sum_{\expand{\KK} \in \K^*} \sum \{\rho_\SS(\{\SS^*\}) : \expand{\SS} \in \LEKstar, \SS^* = (\KK^* : \SS_1^*, \ldots, \SS_a^*)\} \\
	&= \sum \{\rho_\SS(\{\SS^*\}) : \expand{\SS} \in \LEKstar \} = 1,
\end{align*}
where the third equality follows from (\textbf{P}) for $\rho$.

(\textbf{E}) Let $\KK \leq \LL$ be structures in $\K$ and let $\expand{\KK} \in \K^*$. Let $\SS \leq \TT$ be structures in $\LEK$, with $\SS = (\KK : \SS_1, \ldots, \SS_a)$ and $\TT = (\LL : \TT_1, \ldots, \TT_b)$. Let $I := \{1 \leq i \leq b : \SS_i \leq \TT_i\}$ and $J := \{1, \ldots, b\} \setminus I$. Then we have
\begin{align*}
	&  \sum \{\mu_\LL (\{\LL^*\}) : \expand{\KK} \leq \expand{\LL} \in \K^*\} 																							\\
	&= \sum \{\mu_{\LL,\TT} (\{\TT^*\}) : \expand{\KK} \leq \expand{\LL} \in \K^*\} 																				\\
	&= \sum_{\expand{\KK} \leq \expand{\LL} \in \K^*} 
					\sum \{\rho_\TT (\{\TT^*\}) : \expand{\TT} \in \LEKstar, \TT^* = (\LL^* : \TT_1^*, \ldots, \TT_b^*)\}						\\
	&= \sum_{\substack{\expand{\SS} \in \LEKstar \\ \SS^* = (\KK^* : \SS_1^*, \ldots, \SS_a^*)}} 
					\sum \left\{ \rho_\TT (\{\TT^*\}) :
					\begin{aligned}
						&\scriptstyle{\expand{\SS} \leq \expand{\TT} \in \LEKstar }\\
	 					&\scriptstyle{\TT^* = (\LL^* : \TT_1^*, \ldots, \TT_b^*)}
					\end{aligned}
				\right\} \\
	&= \sum_{\substack{\expand{\SS} \in \LEKstar \\ \SS^* = (\KK^* : \SS_1^*, \ldots, \SS_a^*)}} \rho_\SS (\{\SS^*\})
	= \mu_{\KK, \SS} (\{\KK^*\})
	= \mu_\KK (\{\KK^*\}),
\end{align*}
where the fourth equality is by (\textbf{E}) of $\rho$. So we have shown that $\Aut(\Flim(\K))$ is amenable by Proposition \ref{threeone}.
\ep\medbreak

\proc{Proof of (i), $\Rightarrow \Aut(\Flim(\L))$ is amenable.} 

Assume that $\Aut(\Flim(\LEK))$ is amenable. By Proposition (\ref{threeone}) there is a consistent random expansion $\rho$ on $\LEK$. We will show that $\Aut(\Flim(\L))$ is amenable.

Let $\PP$ be a one point structure in $\K$ and let $\expand{\PP} \in \K^*$. For $\LL \in \L$ there is an $\SS \in \LEK$ such that $\SS = (\PP : \LL)$, and every $\expand{\LL}\in\L^*$ gives us an $\SS^* = (\PP^* : \LL^*)$ such that $\expand{\SS} \in \LEKstar$. Using the consistent random expansion $\mu$ we previously defined, we introduce
\[
	\gamma_\LL (\{\LL^*\}) := \rho_\SS (\{\SS^*\}) \cdot \frac{1}{\mu_{\PP} (\{\PP^*\})}.
\]

Note that we must have $\mu_\PP (\{\PP^*\}) \neq 0$ since $\rho_\SS (\{\SS^*\}) \neq 0$ for all $\SS \in \LEK$ and $\expand{\SS} \in \LEKstar$. We prove that $(\gamma_\LL)$ is a consistent random expansion of $\L$ by checking (\textbf{P}) and (\textbf{E}).

(\textbf{P}) Observe that
\begin{align*}
	&  \sum \{\gamma_\LL (\{\LL^*\}) : \expand{\LL}\in\L^*\} 																																													\\
	&= \sum_{\expand{\LL}\in\L^*} \{\rho_\SS (\{\SS^*\}) \cdot \frac{1}{\mu_{\PP} (\{\PP^*\})} : \expand{\SS} \in \LEKstar, \SS^* = (\LL^* : \PP^*)\} 	\\
	&= \frac{1}{\mu_{\PP} (\{\PP^*\})} \cdot \sum_{\expand{\LL}\in\L^*} \{\rho_\SS (\{\SS^*\}) : \expand{\SS} \in \LEKstar, \SS^* = (\PP^* : \LL^*)\} 	\\
	&= \frac{1}{\mu_{\PP} (\{\PP^*\})} \cdot \mu_\PP (\{\PP^*\}) = 1,
\end{align*}
where the third equality follows from the definition of $\mu$.

(\textbf{E}) Let $\LL \leq \KK$ be structures in $\L$. Then there are $\SS \leq \TT \in \LEK$ such that $\SS = (\PP : \LL)$ and $\TT = (\PP : \KK)$. Let $\expand{\LL} \in \L^*$. Then we have the following
\begin{align*}
	&  \sum \{\gamma_\KK (\{\KK^*\}) : \expand{\LL} \leq \expand{\KK} \in \L^*\} 																																	\\
	&= \sum \left \{\frac{1}{\mu_\PP (\{\PP^*\})} \cdot \rho_\TT (\{\TT^*\}) : \expand{\LL} \leq \expand{\KK} \in \L^*, \TT^* = (\PP^*:\KK^*)\right\}	\\
	&= \frac{1}{\mu_\PP (\{\PP^*\})} \cdot \sum \{\rho_\TT (\{\TT^*\}) : \expand{\SS} \leq \expand{\TT}, \SS^* = (\PP^*:\LL^*)\} 											\\
	&= \frac{1}{\mu_\PP (\{\PP^*\})} \cdot \rho_\SS(\{\SS^*\}) 	
			= \frac{1}{\mu_\PP (\{\PP^*\})} \cdot \mu_\PP (\{\PP^*\}) \cdot \gamma_\LL (\{\LL^*\}) = \gamma_\LL (\{\LL^*\}),
\end{align*}
where the third equality follows from (\textbf{E}) for $\rho$. This finishes the verification that $\Aut(\Flim(\LEK))$ is amenable.
\ep\medbreak

\proc{Proof of (ii), $\Rightarrow$.} Assume that $\Aut(\Flim(\LEK))$ is uniquely ergodic. So we have that $\Aut(\Flim(\L))$ and $\Aut(\Flim(\K))$ are amenable, and by Proposition \ref{threeone} there are consistent random expansions $\mu$ and $\gamma$ on $\K$ and $\L$ respectively. Suppose that one of $\Aut(\Flim(\L))$ or $\Aut(\Flim(\K))$ is not uniquely ergodic. Then there is a consistent random expansion $\mu^\prime$ on $\K$ such that $\mu \neq \mu^\prime$ or there is a consistent random expansion $\gamma^\prime$ on $\L$ such that $\gamma \neq \gamma^\prime$. Then there is a structure $\KK \in \K$ and an expansion $\expand{\KK}\in\K^*$ such that
\[
	\mu_\KK (\{\KK^*\}) \neq \mu_\KK^\prime(\{\KK^*\})
\]
or there is a structure $\LL \in \L$ and an expansion $\expand{\LL}\in\L^*$ such that
\[
	\gamma_\LL(\{\LL^*\}) \neq \gamma_\LL^\prime(\{\LL^*\}).
\]
Now consider the structure $\SS = (\KK : \LL, \ldots, \LL) \in \LEK$ with expansion $\expand{\SS}\in\LEKstar$ where $\SS^* = (\KK^* : \LL^*, \ldots, \LL^*)$, with $a := \vert \KK \vert$ many $\LL$. Using similar arguments to the proof of $[(1), \Leftarrow]$ of this theorem, we have that $\gamma\otimes\mu, \gamma\otimes\mu^\prime$ and $\gamma^\prime \otimes\mu$ are consistent random expansions on $\LEK$. In particular if $\mu \neq \mu^\prime$, we have
\begin{align*}
	(\gamma\otimes\mu)_\SS (\{\SS^*\}) 	&= \mu_\KK (\{\KK^*\}) \cdot \prod_{i=1}^a \gamma_{\LL} (\{\LL^*\}) \\
																	&\neq \mu^\prime_\KK (\{\KK^*\}) \cdot \prod_{i=1}^a \gamma_{\LL} (\{\LL^*\}) 
																	= (\gamma\otimes \mu^\prime)_\SS (\{\SS^*\}),
\end{align*}
and if $\gamma \neq \gamma^\prime$, then we have
\begin{align*}
	(\gamma\otimes\mu)_\SS (\{\SS^*\}) 	&= \mu_\KK (\{\KK^*\}) \cdot \prod_{i=1}^a \gamma_{\LL} (\{\LL^*\}) \\
																			&\neq \mu_\KK (\{\KK^*\}) \cdot \prod_{i=1}^a \gamma^\prime_{\LL} (\{\LL^*\}) 
																			= (\gamma^\prime \otimes\mu)_\SS (\{\SS^*\}).
\end{align*}
Therefore we have two distinct consistent random expansions on $\LEK$. This is in contradiction to the unique ergodicity of $\LEK$, according to Theorem~\ref{threeone}, so $\Aut(\Flim(\L))$ and $\Aut(\Flim(\K))$ must be uniquely ergodic.
\ep\medbreak

\proc{Proof of (ii), $\Leftarrow$.} Now assume that $\Aut(\Flim(\L))$ and $\Aut(\Flim(\K))$ are uniquely ergodic, and let $\mu$ and $\gamma$ be the unique consistent random expansions on $\K$ and $\L$ respectively. According to the first part of this theorem there is a consistent random expansion on $\LEK$. We will show that $\rho := \gamma\otimes \mu$ is the unique consistent random expansion on $\LEK$, as defined in the previous part of the proof.

Let $\SS = (\KK : \SS_1, \ldots, \SS_a)$ be a structure in $\LEK$ with expansion $\expand{\SS} \in \LEKstar$ given by $\SS^* = (\KK^*:\SS_1^*, \ldots, \SS_a^*)$. In the previous proof of [(1), $\Rightarrow \Aut(\Flim(\K))$ is amenable] we described a consistent random expansion on $\K$ given by $\rho$, so we may notice that by unique ergodicity, this is exactly $\mu$. 

Fix any $\TT_2, \ldots, \TT_a \in \L$ which will be used to define measures on $\L$. They can be $\SS_2, \ldots, \SS_a$ if you like, but for purposes of clarity we use $\TT_2, \ldots, \TT_a$.

Define
\[
	p_0 := \mu_\KK (\{\KK^*\}).
\]

Let $\LL \in \L$ with $\expand{\LL}\in\L^*$ be given. Consider the map
\begin{align*}
	\gamma_{\LL}^1 (\{\LL^*\}) &:= \frac{1}{p_0} \cdot
		\sum \left\{\rho_{\XX_1} (\XX_1^*) :
						\begin{aligned}
							&\scriptstyle{\XX_1 = (\KK : \LL, \TT_2, \ldots, \TT_a) \in \LEK}\\
							&\scriptstyle{\XX_1^* = (\KK^*:\LL^*, \TT_2^*, \ldots, \TT_a^*) \in \LEKstar}
						\end{aligned}
					\right\} 
\end{align*}

Notice that the sum does not run over $\LL^*$ and $\KK^*$, which are fixed.

In a similar way as the proof of [(1), $\Rightarrow \Aut(\Flim(\L))$ is amenable] we may conclude that $\gamma^1$ is a consistent random expansion on $\L$ that does not depend on the choice of $\TT_2, \ldots, \TT_a$. Since $\Aut(\Flim(\L))$ is uniquely ergodic we must have that $\gamma^1 = \gamma$. In particular, for $\LL = \SS_1$ we can define
\[
	p_1 := p_0 \cdot \gamma_{\SS_1} (\{\SS_1^*\}).
\]

Now consider the map $\gamma^2$, for $\LL \in \L$ with $\expand{\LL} \in \L^*$ given by:
\begin{align*}
	\gamma_{\LL}^2 (\{\LL^*\}) &:= \frac{1}{p_1} \cdot
		\sum \left\{\rho_{\XX_2} (\{\XX_2^*\}) :
						\begin{aligned}
							&\scriptstyle{\XX_2 = (\KK:\SS_1, \LL, \TT_3, \ldots, \TT_a) \in \LEK}\\
							&\scriptstyle{\XX_2^* = (\KK^*:\SS_1^*, \LL^*, \TT_3^*, \ldots, \TT_a^*) \in \LEKstar}
						\end{aligned}
					\right\} 
\end{align*}

Notice that the sum does not run over $\SS_1^*, \LL^*$ and $\KK^*$, which are fixed.

Again similar arguments as in the previous proof show that $\gamma^2$ is a consistent random expansion on $\L$, and by unique ergodicity we have that $\gamma_2 = \gamma$. In particular we define:
\[
	p_2 := p_1 \cdot \gamma_{\SS_2} (\{\SS_2^*\}).
\]

Continuing on in this way we obtain
\[
p_a := \mu_\KK (\{\KK^*\}) \cdot \gamma_{\SS_1}(\{\SS_1^*\}) \cdot \ldots \cdot \gamma_{\SS_a}(\{\SS_a^*\})
\]
with $p_a = \rho_\SS (\{\SS^*\})$. Therefore we have proved that $\rho = \gamma \otimes \mu$, so it must be unique.
\ep\medbreak

%
\subsection{The corollaries.}

For $n \leq \omega$, define $[n] := \{i \in \omega : i < n\}$.

Let $\TT$ be one of the tournaments $\QQ, \SS(2), \TT^\omega$ or $\CC_3,$ and let $T$ be the underlying set of $\TT$. For $n \leq \omega$ we denote by $\TT[\II_n]$ the directed graph with the underlying set $T \times [n]$ and the edge relation given by
\[
	(x,i) \rightarrow (y,j) \text{ iff } x \rightarrow y
\]
and $\II_n[\TT]$ the tournament with the underlying set $[n]\times T$ and edge relation given by
\[
	(i,x) \rightarrow (j,y) \text{ iff } (i=j, x \rightarrow y)
\]

Consider $\II_n$ as a structure on the empty signature. Therefore we have the following:

\begin{cor}\label{cor:LEK_cor} Let $\TT$ be one of the tournaments $\QQ, \SS(2), \TT^\omega$ or $\CC_3$, and let $n \leq \omega$. Then,
	\begin{enumerate}
		\item $\Aut(\SS(2)[\II_n])$ is not amenable. For $\TT \neq \SS(2)$, $\Aut(\TT[\II_n])$ is uniquely ergodic.
		\item $\Aut(\II_n [\SS(2)])$ is not amenable. For $\TT \neq \SS(2)$, $\Aut(\II_n [\TT])$ is uniquely ergodic.
	\end{enumerate}
\end{cor}

\proc{Proof.} Non-amenability of $\Aut(\SS(2))$ comes from unpublished communication with A.Kechris, and the argument is given in Theorem \ref{thm:Sn_amen}.

In the case that $\TT[\II_n]$ or $\II_n[\TT]$ is a finite structure then its automorphism group is finite and therefore uniquely ergodic. So let us assume that they are infinite. Thus we may view the structures as Fra\"iss\'e limits of the form $\LLEKK := \Flim(\LEK)$ where both $\L$ and $\K$ are Fra\"iss\'e structures, or one of them (e.g. $\II_n$) is not a Fra\"iss\'e structure simply because it is not an infinite structure. This can be rectified by going through the proof of Theorem~\ref{A_LEK} and noticing that the assumption that $\L$ and $\K$ are infinite is not used.
\ep\medbreak

\section{Ramsey expansions of countable homogeneous directed graphs}

We now introduce the directed graphs that appear in Cherlin's classification. More detailed descriptions will be given in the relevant sections.

%
\subsection{Summary and known results.}
Denote by $\II_n$ the edgeless directed graph on $n$ vertices, where $n \leq \omega$.

Denote by $\CC_3$ the directed $3$-cycle. Specifically, $\CC_3 = (C_3, \rightarrow^{C_3})$ is the directed graph such that $C_3 = \{a,b,c\}$ with $a \rightarrow^{C_3} b, b \rightarrow^{C_3} c$ and $c \rightarrow^{C_3} a$.

The following is the classification of countable homogeneous directed graphs (see \cite{Cher98}). The infinite structures here are Fra\"iss\'e structures.

\begin{enumerate}
	\item The finite digraphs $\CC_3$ and $\II_n$ for $n < \omega$.
	\item $\II_\omega$ is the edgeless directed graph on $\omega$ vertices.
	\item $\QQ, \SS(2)$ and $\TT^\omega$ are tournaments. 
		\begin{enumerate}
			\item $\QQ$ is the set of rational numbers where we take $x \rightarrow^\QQ y$ iff $x < y$.
			\item $\SS(2)$ is the \textbf{dense local ordering} which may be seen as the set of points on the unit circle with rational arguments such that $e^{i \theta} \rightarrow^{S(2)} e^{i \phi}$ iff $0 < \phi - \theta < \pi$. See Section~\ref{sec:def_S(n)}.
			\item $\TT^\omega$ is the generic tournament, i.e. the Fra\"iss\'e limit of the class of all finite tournaments.
		\end{enumerate}
	\item $\TT[\II_n], \II_n[\TT]$, where $n \leq \omega$, and $\TT$ is one of the tournaments in (3) or $\CC_3$. This is a type of ``$n$-point cover (or blowup)" of the nodes of the tournament $\TT$. See Section~\ref{sec:LEK}.
	\item $\hat{\TT}$, for $\TT = \II_1, \CC_3, \QQ$ or $\TT^\omega$, is a type of ``two point cover (or blowup)" of the vertices of a tournament. See Section~\ref{sec:def_blowups}.
	\item $\DD_n$, for $1 < n \leq \omega$, the complete $n$-partite directed graph with countably many nodes. See Section~\ref{sec:def_Dn}.
	\item $\SS$ is the semi-generic graph. See Section~\ref{sec:def_semigeneric}.
	\item $\SS(3)$ is a directed graph which may be seen as the set of points on the unit circle with rational arguments such that $e^{i \theta} \rightarrow^{S(3)} e^{i \phi}$ iff $0 < \phi - \theta < \frac{2 \pi}{3}$. This is not a tournament. See Section~\ref{sec:def_S(n)}.
	\item $\PP$ is the generic poset, i.e. the Fra\"iss\'e limit of the class of all finite posets such that $x \rightarrow^P y$ iff $x <^P y$.
	\item $\PP(3)$ is the ``twisted'' generic poset in three parts. See Section~\ref{sec:def_P(3)}.
	\item $\GG_n$, for $n > 1$, is the generic directed graph with the property that $\II_{n+1}$ can't be embedded in $\GG_n$. See Section~\ref{sec:hypergraph}.
	\item $\FF(\T)$ is the Fra\"iss\'e limit of the class of finite directed graphs which do not embed any member of $\T$, where $\T$ is a fixed set of finite tournaments each of which has at least three vertices; they are \textbf{F}orbidden. See Section~\ref{sec:hypergraph}.
\end{enumerate}

In this paper we go through Cherlin's classification and examine the automorphism group of each of these structures with respect to amenability and unique ergodicity. We consider each automorphism group as a topological group with the pointwise convergence topology, see \cite{BK96} for more details.

First we give a list of known facts:

\begin{enumerate}
	\item $\Aut(\II_\omega) \cong S_\infty$ is amenable and uniquely ergodic, see \cite{GW02}, \cite[Proposition 10.1]{AKL12}.
	\item $\Aut(\QQ)$ and $\Aut(\TT^\omega)$ are amenable and uniquely ergodic, see \cite{pes98}, \cite[Theorem 6.1, Theorem 2.2]{AKL12}.
	\item $\Aut(\SS(2))$ and $\Aut(\SS(3))$ are not amenable, see \cite[Theorem 3.1]{ZUC13} and private communication with Kechris.
	\item $\Aut(\PP)$ is not amenable, see \cite[Section 3]{KSOK12}.
\end{enumerate}

In addition, the essential arguments for the amenability and unique ergodicity of $\Aut(\GG_n)$ and $\Aut(\FF(\T))$ are contained in \cite[Theorem 5.1]{AKL12} which shows a similar theorem for undirected graphs (with some combinatorial conditions). Amenability and unique ergodicity of $\Aut(\GG_n)$ and $\Aut(\FF(\T))$ are not direct corollaries of their theorem or their proof, modifications had to be made.

%
\subsection{Weak local orders, $\SS(2), \SS(3)$ and $\SS(n)$.}
\label{sec:def_S(n)}

Fix $2 \leq n < \omega$. The directed graph $\SS(n)$ may be seen as the set of points on the unit circle with rational arguments such that $e^{i \theta} \rightarrow^{S(n)} e^{i \phi}$ iff $0 < \phi - \theta < \frac{2 \pi}{n}$. This is a tournament iff $n=2$, and is homogeneous iff $n=2,3$.

Let $I_k$, for $0 \leq k \leq n-1$ be unary relational symbols. We consider the structure $\SS(n)^*$ in the signature $\{\rightarrow, I_0, \ldots, I_{n-1}\}$ where:
	\begin{itemize}
		\item $\SS(n)^* \vert \{\rightarrow\} = \SS(n)$,
		\item $I_k^{\SS(n)}(x) \Leftrightarrow x \in e^{i\theta}$ and $\theta \in (\frac{k \cdot 2\pi}{n}, \frac{(k+1) \cdot 2\pi}{n})$,
	\end{itemize}

Set $\mathcal{S}(n) := \Age(\SS(n))$ and $\mathcal{S}(n)^* := \Age(\SS(n)^*)$.

For $n \geq 4$, $\SS(n)$ is not a Fra\"iss\'e structure. For $n=2,3$, $(\mathcal{S}(n), \mathcal{S}(n)^*)$ is an excellent pair, see \cite{LNS10,NVT13}.

%
\subsection{Partial order $\PP$ and a variant $\PP(3)$.}
\label{sec:def_P(3)}

Let $\P$ be the class of finite posets in signature $\{\leq\}$. Let $P_0, P_1, P_2$ be unary relational symbols and let $\P_3$ be the class of finite structures of the form $(A, \leq^A, P_0^A, P_1^A, P_2^A )$ where $(A, \leq^A) \in \P$ and $A = \bigsqcup_{i=0}^2 \{x : P_i^A(x)\}$.

It is easy to see that $\P_3$ is a Fra\"iss\'e class with limit $\PP_3 = (R, \leq^R, P_0^R, P_1^R, P_2^R)$. Using $\PP_3$ we define the structure $\PP(3) = (R, \rightarrow^R)$ such that for $x,y \in R$ with $P_i^R(x)$ and $P_j^R(y)$ we have $x \rightarrow^R y$ iff one of the following conditions is satisfied:

\begin{enumerate}
	\item $j = i$ and $x <^R y$; or
	\item $j = i+1 \mod 3$ and $y <^R x$; or
	\item $j = i+2 \mod 3$ and $x$ is $<^R$-incompatible with $y$.
\end{enumerate}

\begin{figure}[!ht]
    \centering
    \includegraphics[width=0.4\textwidth]{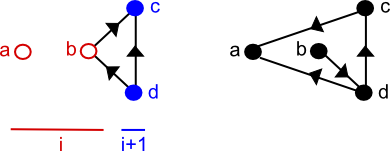}
    \caption{A structure in $\P_3$ and its corresponding structure in $\P(3)$.}
\end{figure}

The structure $\PP(3)$ is a Fra\"iss\'e structure, see \cite{Cher98}, which is called the generic twisted poset. Note that each $P_i$ induces a copy of $\PP$ that is cofinal in $\PP(3)$.  We also consider $\PP(3)^* = (R, \rightarrow^R, P_0^R, P_1^R, P_2^R, \preceq^R)$, which is also a Fra\"iss\'e structure, where $\preceq^R$ is a linear order on $R$ that extends the partial order $(R, \leq^R)$ given by untwisting $(R, \rightarrow^R, P_0^R, P_1^R, P_2^R)$. 

In what follows we will use $\leq$ for untwisted partial orders, $\rightarrow$ for the corresponding twisted directed graph, and $\preceq$ for the linear order that extends $\leq$. We will not refer to an untwisted partial order's natural directed graph.

Let $\P(3) := \Age(\PP(3))$ and $\P(3)^* := \Age(\PP(3)^*)$. The pair $(\P(3), \P(3)^*)$ is an excellent pair, see \cite{JLNW14}.

%
\subsection{2-covers of tournaments, $\hat{\QQ}$ and $\hat{\TT^\omega}$.}
\label{sec:def_blowups}

Here we discuss a way of 2-covering (or blowing up) points of a tournament so that it has much of the same structure, but it is no longer a tournament. In the case of covering $\TT^\omega$, the Ramsey expansion is straightforward, being essentially convex linear orders. In the case of covering $\QQ$, the Ramsey expansion is more subtle and must interact suitably with the linear order of $\QQ$. We explain these expansions in some depth because their discussion was unexpectedly absent from \cite{JLNW14}; in particular we show that these expansions are indeed Ramsey expansions.

\subsubsection{The 2-cover $\hat{\TT}$.}

Let $\TT$ be one of the following tournaments: $\II_1, \CC_3, \QQ$ or $\TT^\omega$, and let $T$ be its underlying set, and let $\T := \Age(\TT)$. We consider the structure $\hat{\TT}$ with underlying set\footnote{\textbf{C}astor and \textbf{P}ollux are the Gemini twins. Using the set $\{C,P\}$ should help the reader, in the proofs that follow, distinguish what notation represents variables and what represents fixed objects} $\hat{T} = T \times \{C,P\}$ and edge relation $\rightarrow^{\hat{T}}$ given by:
\[
	(x,i) \rightarrow^{\hat{T}} (y,j) \Leftrightarrow \left((x \rightarrow^T y, i\neq j) \vee (y \rightarrow^T x, i=j)\right).
\]

\begin{figure}[!ht]
    \centering
    \includegraphics[scale=0.5]{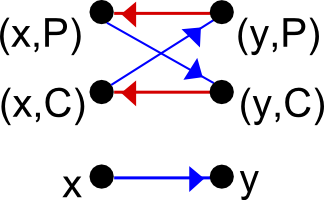}
    \caption{The bottom graph is an example of $T \in \T$, the upper graph is $\hat{T} \in \hat{\T}$.}
\end{figure}

\subsubsection{Description of $\hat{\T}$ and $\T^*$.}

For $\TT \in \{\II_1, \CC_3\}$ the structure $\hat{\TT}$ is finite so $\Aut(\hat{\TT})$ is finite, hence uniquely ergodic. For $\TT \in \{\QQ, \TT^\omega\}$, $\hat{\TT}$ is a Fra\"iss\'e structure with corresponding Fra\"iss\'e class $\hat{\T} := \Age(\hat{\TT})$. 

Let $\T^*$ be the collection of structures of the form $(A, \rightarrow^A, \leq^A)$ where $(A, \rightarrow^A) \in \T$ and $\leq^A$ is a linear order on $A$. So $(\T, \T^*)$ is an excellent pair, by \cite{AH78} and \cite{NR77,NR83,NR89}.

\subsubsection{Description of $\hat{\T}^*$.}

For each structure $(A, \rightarrow^A) \in \hat{\T}$, the relation $\perp^A$ is an equivalence relation which gives the partition $A = A_1 \sqcup \ldots \sqcup A_k$ where each class has at most two elements. In the following we describe a Fra\"iss\'e class $\hat{\T}^*$ such that $(\hat{\T},\hat{\T}^*)$ is an excellent pair.

Let $\hat{\T}^*$ contain structures of the form $(A, \rightarrow^A, \leq^A, I_0^A, I_1^A)$ where:
\begin{itemize}
	\item $(A, \rightarrow^A) \in \hat{\T}$, 
	\item $\leq^A$ is a linear order on $A$,
	\item $I_0^A$ and $I_1^A$ are unary relations on $A$ which partition $A$, and
	\item each $A_i$ is an interval with respect to $\leq^A$.
\end{itemize}

Though a slight abuse of notation, we denote by $\leq^A$ the linear ordering induced by $\leq^A$ on the set $\{A_1, \ldots, A_k\}$. The correct expansion of $\hat{\T}$ for $\T = \Age(\TT^\omega)$ allows arbitrary partitions by $I_0$ and $I_1$ in each $\perp$-equivalence class, so long as the smallest element in each equivalence class, with respect to the linear ordering, belongs to $I_0$. For $\T = \Age(\QQ)$ there is an additional subtlety which we will explain in the following section; essentially the expansion must be given by a transversal that coheres with the inherent linear order of $\QQ$, but for $\TT^\omega$ there is no such linear order to cohere with. This is similar, although not identical, to how the correct expansions for the generic partial orders are the linear orders that extend the partial orders. 

The proofs of the following two facts are relegated to the appendix.

\begin{theo}[\ref{thm:hatT_Ramsey}] $\hat{\T}^*$ is a Ramsey class.
\end{theo}

\begin{prop}[\ref{prop:hatT_EP}] $\hat{\T}^*$ satisfies the \textbf{EP} with respect to $\hat{\T}$.
\end{prop}

\subsubsection{Description of $\hat{\T}^*$ for $\T = \Age(\QQ)$.}

Let $\AA = (A, \rightarrow^A) \in \hat{\T}$ be a structure with $k$ many $\perp^A$-equivalence classes $A_1, \ldots, A_k$ such that $\vert A_i \vert =2$ for all $i \leq k$. Suppose that $A_i = \{(i,C), (i,P)\}$ and that for $i \neq j$ we have:
\begin{itemize}
	\item $(i,C) \rightarrow^A (j,C) \Leftrightarrow i > j$;
	\item $(i,P) \rightarrow^A (j,P) \Leftrightarrow i > j$;
	\item $(i,P) \rightarrow^A (j,C) \Leftrightarrow i < j$;
	\item $(i,C) \rightarrow^A (j,P) \Leftrightarrow i < j$.
\end{itemize}

\begin{figure}[!ht]
    \centering
    \includegraphics[scale=0.5]{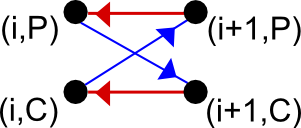}
    \caption{Two columns $A_i$ and $A_{i+1}$ of $A$.}
\end{figure}

Now we examine when a transversal $T \sse \AA$ forms a linear ordering, that is, a linear order on $T$ that also gives rise to the induced subgraph on $T$. The following lemma says that $T$ forms a linear order so long as there is at most one change of levels, and no zigzags. It also establishes that if a structure $\AA$ has $k$ equivalence classes each with two points, then $\AA$ has exactly $2k$ expansions.

\begin{figure}[!ht]
		\centering
    \includegraphics[scale=0.5]{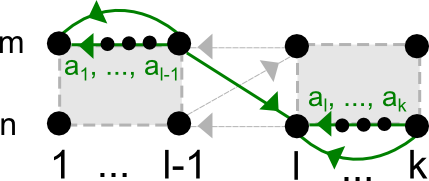}
		\caption{A linear order where $a_{l-1} < \ldots < a_1 < a_k < \ldots < a_l$.}
		\label{fig:l}
\end{figure}

\begin{lemma}\label{lem:l}A sequence of vertices $(a_i)_{i=1}^k$ with $a_i \in A_i$ forms a linear ordering iff
\begin{enumerate}
	\item All $a_i$ have the same second coordinate; or
	\item There is $l < k$ such that for all $i < l$, the $a_i$ have the same second coordinate $m$, and for all $i \geq l$ the $a_i$ have the same second coordinate $n \neq m$. (See Figure~\ref{fig:l}.)
\end{enumerate}
\end{lemma}

\proc{Proof.} It is enough to consider the following. Let $i < j < k$ and $m \neq n$. Then we have:
\[
	(i,m) \rightarrow^A (j,n) \rightarrow^A (k,m) \rightarrow^A (i,m)
\]
so the directed graph induced by $(i,m), (j,n), (k,m)$ is not a linear ordering. 

\begin{figure}[!ht]
    \centering
    \includegraphics[scale=0.5]{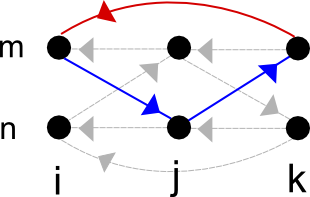}
    \caption{The following edges have been omitted for readability: $(i,m)$ to $(k,n)$ and $(i,n)$ to $(k,m)$.}
\end{figure}

\ep\medbreak

We define an expansion class $\hat{\T}^*$ using such a sequence $\vec{a} = (a_i)_{i=1}^k$ in $A$, which forms a linear ordering $\leq^{\vec{a}}$. Then we introduce indicators $I_0^A$ and $I_1^A$ such that for $x \in A$ we have:
\begin{align*}
	I_1^A (x) &\Leftrightarrow x \in \{a_1, \ldots, a_k\} \\
	I_0^A (x) &\Leftrightarrow x \notin \{a_1, \ldots, a_k\} \\
\end{align*}

We define a linear ordering $\leq^A$ on $A$ such that for $x \in A_i$ and $y \in A_j$ we have:
\[
	x <^A y \Leftrightarrow \left((i=j,I_1^A(x)) \vee (a_i <^{\vec{a}} a_j)\right).
\]

In the case where some columns of $\AA$ do not have two elements, we are a little more careful. For every structure $\AA \in \hat{\T}$ there is a unique structure $\BB \in \hat{\T}$, up to isomorphism, which contains the same number of $\perp^B$-equivalence classes of $\AA$ each of which has exactly two elements. Thus we may define an expansion of $\BB$, then by taking the restriction to $\AA$, we get an expansion of $\AA$.

\begin{lemma}\label{lem:hatT_iso} For all $\AA \in \hat{\T}$ and all $\expand[a]{\AA}, \expand[b]{\AA} \in \hat{\T}^*$, the structures $\expand[0]{\AA}$ and $\expand[1]{\AA}$ are isomorphic. Therefore $(\hat{\T}, \hat{\T}^*)$ satisfies the \textbf{EP}.
\end{lemma}

\proc{Proof.} Let $\vec{b} = (b_i)_{i=1}^k$ be a sequence in $A$ given by $b_i = (i,C)$ for each $i \leq k$. Then $\vec{b}$ forms a linear ordering and it induces an expansion of $\AA$, call it $\AA^b = (A, I_0^b, I_1^b, \leq^b)$ (where the vector notation is dropped for readability). Let $\vec{a}$ be another sequence in $A$ which induces the expansion $\AA^a$.

Consider the map $\pi_{\vec{a}} : A \longrightarrow A$ given by:
\[
	\pi_{\vec{a}}(i,m) = 
		\begin{cases}
			(i+l, m)  &: 1   \leq i \leq l \\
			(i-l,1-m) &: l+1 \leq i \leq k 
		\end{cases}
\]
where $l$ is given by Lemma~\ref{lem:l}.

It is easy to see that $\pi_{\vec{a}}$ is an automorphism of $\AA$ and moreover that $\pi_{\vec{a}}$ is an isomorphism between $\expand[a]{\AA}$ and $\expand[b]{\AA}$. Therefore any two expansions of $\AA$ are isomorphic.
\ep\medbreak

%
\subsection{Complete $n$-partite directed graph, $\DD_n$.}
\label{sec:def_Dn}

For $n\in\NN$ let $\D_n$ be the class of finite digraphs $(A, \rightarrow^A)$ in which $\perp^A$ is an equivalence relation with at most $n$ many equivalence classes. We will also consider $\D_\omega := \bigcup_{n < \omega} \D_n$. In this way we obtain Fra\"iss\'e classes $\D_n$, for $n \leq \omega$, with corresponding Fra\"iss\'e limits $\DD_n$, for $n \leq \omega$. 

We denote by $\D_\omega^*$ the class of finite structures of the form $(A, \rightarrow^A, \leq^A)$ where $(A, \rightarrow^A) \in \D_\omega$ and $\leq^A$ is a linear order on $A$ such that $\forall x,y,z \in A$ we have
\[
	x <^A y <^A z, x \perp^A z \Rightarrow x \perp^A y \perp^A z,
\]
which is a type of \textbf{convexity}.

For a finite $n$, we let $\D_n^*$ be the class of finite structures of the form $(A, \rightarrow^A, \leq^A, \{I_i^A\}_{i=1}^n)$ where $(A, \rightarrow^A, \leq^A) \in \D_\omega^*$ and each $I_i^A$ is a unary relation on $A$ such that $\forall x,y \in A$ we have:
\begin{itemize}
	\item $(\exists i \leq n)(I_i^A(x))$,
	\item $I_i^A(x), x \perp^A y \Rightarrow I_i^A(y)$,
	\item $I_i^A(x), I_j^A(y), i < j \Rightarrow x<^A y$.
\end{itemize}

The $I_i^A$ indicate the $n$-parts, and the ordering is convex with respect to the parts.

%
\subsection{Semi-generic multipartite digraph $\SS$ and variants $\SS_\leq, \SS_R$.}
\label{sec:def_semigeneric}

\subsubsection{Description of $\SS$.}

Let $\mathcal{S}$ be the class of finite directed graphs of the form $(S, \rightarrow^S)$ with the following properties:
	\begin{enumerate}
		\item The binary relation $\perp^S$ defined on $S$ by $x \perp^S y \Leftrightarrow \neg(x \rightarrow^S y \vee y \rightarrow^S x)$ is an equivalence relation on $S$. (We call the equivalence classes \textbf{columns}.)
		\item For $x,y,t_x,t_y \in S$, where $x \perp^S t_x$ and $y \perp^S t_y$, we have that the number of edges directed from $\{x,t_x\}$ to $\{y,t_y\}$ is even.
	\end{enumerate}

\begin{figure}[!ht]
    \centering
    \includegraphics[width=0.8\textwidth]{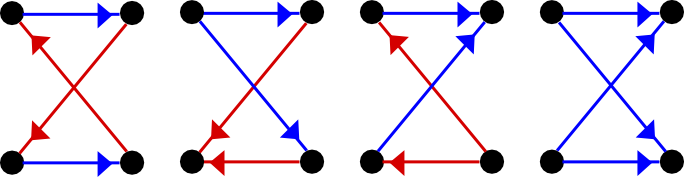}
    \caption{The 4 possible digraphs, with 2 nodes on two columns, with the parity condition, up to reflection. The first is in ``general position''.}
\end{figure}

Condition (1) ensures that the digraphs are complete $n$-partite, for some $n$. The parity condition (2) might seem artificial, but it has the following nice property which says ``If you know three edges, then you know the fourth edge''.

\begin{lemma}[Three of four] Let $\AA = (\{x,y,t_x,t_y\}, \rightarrow^S)$, with $x \perp^S t_x$, $y \perp^S t_y$, $\neg(t_x \perp^S t_y)$. If 3 of the directed edges between $\{x,t_x\}$ and $\{y,t_y\}$ are specified, then there is a unique directed edge between $\{x,t_x\}$ and $\{y,t_y\}$ that, when added, satisfies the parity condition.
\end{lemma}

\begin{figure}[!ht]
    \centering
    \includegraphics[scale=0.5]{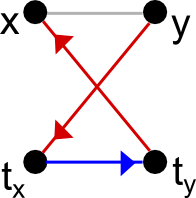}
    \caption{An example of the three of four lemma.}
\end{figure}
	
Let $\mathcal{S}^*$ be the class of finite structures of the form $(A, \rightarrow^A, R^A, \leq^A)$ where $(A, \rightarrow^A) \in \mathcal{S}$, $R^A$ is a binary relation on $A$ and $\leq^A$ is a linear ordering on $A$ with the property that:
	\begin{enumerate}
		\item If $a$ is the number of $\perp^A$-equivalence classes, then there is a linear ordering $T = \{t_1 \leq t_2 \leq \ldots \leq t_a\}$, called a \textbf{transversal}, which we consider as a directed graph $\TT = (T, \rightarrow^T) \in \mathcal{S}$ given by $t_i \rightarrow t_j \Leftrightarrow t_i < t_j$. Then there is a $\BB = (B, \rightarrow^B) \in \mathcal{S}$ such that $\AA \leq \BB, \TT \leq \BB$ and $\BB$ also has $a$ many $\perp^B$-equivalence classes. Also, $\TT$ must be defined on each of the columns of $\AA$ (that is, $\forall x \in \AA, \exists t_i$ such that $x \perp t_i$). See Lemma \ref{lem:semi_amalgam} for further discussion.
		\item If $R^A(x,y)$, then $\neg(x \perp^A y)$. If $t_i \perp^B x$, then we have $$R^A (x,y) \Leftrightarrow t_i \rightarrow^B y.$$
		\item If $x \perp^A z$ and $x <^A y <^A z$ then $x \perp^A y \perp^A z$. If $t_i \perp^B x$ and $t_j \perp^B y$, then $x <^A y \Leftrightarrow t_i < t_j$. (This is a type of \textbf{convexity}.)
	\end{enumerate}

The condition (2), and the three of four lemma, ensures that the digraph structure of $\TT$ amalgamated with $\AA$ can be reconstructed from the transversal and the relations $R^A(x,y)$ (and \textit{vice versa}). The condition (3) says that the linear order is \textbf{convex} with respect to the $\perp^A$-equivalence classes.

\begin{figure}[!ht]
    \centering
    \includegraphics[scale=0.7]{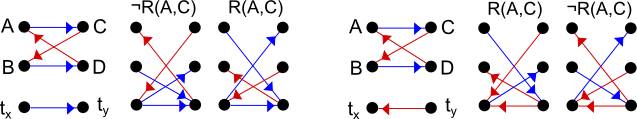}
    \caption{The two possible orientations of a transversal on two columns, and each of their two possible amalgamations with a graph on two columns. The edges between $\{A,B\}$ and $\{C,D\}$ are omitted for readability.}
\end{figure}

The classes $\mathcal{S}$ and $\mathcal{S}^*$ are Fra\"iss\'e classes with limits $\SS$ and $\SS^*$ respectively, see \cite{Cher98,JLNW14}.

\begin{lemma}\label{lem:semi_amalgam} Every structure $\AA \in \mathcal{S}$ can be amalgamated with a transversal with the same columns. Moreover, the transversal can be chosen to respect an arbitrary linear order of the columns of $\AA$.
\end{lemma}

\proc{Proof.} Let $C_i$ for $i \leq a$ be an enumeration of the columns of $\AA$. From each column choose a $c_i \in C_i$. We will amalgamate a transversal $\TT = (T, \rightarrow^T) \in \mathcal{S}$, where $T = \{t_1 \leq t_2 \leq \ldots \leq t_a\}$, and $t_i \perp c_i$ for each $i\leq a$.

First we describe a digraph structure on the nodes $X = \{c_i : i \leq a\} \cup \{t_i : i \leq a\}$. Between $c_i$ and $c_j$ maintain the same edge direction as in $\AA$, and add an edge from $t_i$ to $t_j$ iff $i < j$. For each $i \neq j$ there are many possible choices for the edges between $\{c_i, t_i\}$ and $\{c_j, t_j\}$ so that it respects the parity condition. Note that, in terms of the parity condition, the edges between $\{c_i, t_i\}$ and $\{c_j, t_j\}$ don't interact with the edges to any other columns. Denote by $\XX$ this digraph structure on $X$.

By the Strong Amalgamation Property for $\mathcal{S}$, $\AA$ and $\XX$ can be amalgamated along $\AA \restrict \{c_i : i \leq a\}$, which yields the desired result.
\ep\medbreak

\subsubsection{The Relation $R$.}

Fixing a point $x \in \SS$, the relation $R_x$ induces an equivalence relation with two classes on each other column $A$ which does not contain $x$, where
\[
	R_x (y) \Leftrightarrow R(x,y).
\]

The three of four lemma ensures that if $x \perp x^\prime$, then $R_x$ and $R_{x^\prime}$ induce the same partition on $A$. Thus we refer to the \textbf{partition of a column} $A_i$ given by another column $A_j$. In fact, any $N$ columns $A_1, \ldots, A_N$ induce an equivalence relation on each other column $A$, which contains $2^N$ equivalence classes (some of which may be empty). These column partitions are equivalence relations that are finer than the $\perp$-equivalence relation, but to avoid confusion we shall refer to \textit{column partitions} and $\perp$-\textit{equivalence relations}.

\begin{figure}[!ht]
    \centering
    \includegraphics[scale=0.7]{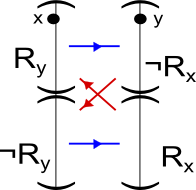}
    \caption{The column partitions two columns induce on each other. Notice the relative positions of $R_x$ and $R_y$.}
		\label{fig:semi_R}
\end{figure}

For more discussion about the relation $R$ and column partitions see \cite[Section 10]{JLNW14}.

\subsubsection{Variants.}
Now we provide some related expansions of $\mathcal{S} := \Age(\SS)$. These variants will allow us to provide partial results towards the unique ergodicity of $\Aut(\SS)$. We have a good understanding of each of the variants separately, but when combined they form $\SS$, which we do not fully understand.

Let $\mathcal{S}_R := \mathcal{S}^* \vert \{\rightarrow, R\}$ and $\mathcal{S}_\leq := \mathcal{S}^* \vert \{\rightarrow, \leq\}$. It is not hard to see that $\mathcal{S}_R$ and $\mathcal{S}_\leq$ are Fra\"iss\'e classes, and that $(\mathcal{S}_R, \mathcal{S}^*)$ and $(\mathcal{S}_\leq, \mathcal{S}^*)$ are excellent pairs. Let $\SS_R := \Flim(\mathcal{S}_R)$ and $\SS_\leq := \Flim(\mathcal{S}_\leq)$.

%
\subsubsection{Generic omitters.}
\label{sec:def_omitters}

Here we introduce two Fra\"iss\'e structures, $\GG_n$ and $\FF(\T)$, that are generic subject to the condition that they omit a specific class of finite directed graphs.

Let $\G_n$ ($n \geq 2$) be the class of finite directed graphs $\GG$ with the property that $\II_{n+1} \not\into \GG$. This is a Fra{\"i}ss{\'e} class with limit $\GG_n := \Flim(\G_n)$. Let $\G_n^*$ be the class of structures of the form $(A, \rightarrow^A, \leq^A)$ with the property that $(A, \rightarrow^A) \in \G_n$ and $\leq^A$ is a linear order on $A$. Then $\G_n^*$ is a Fra{\"i}ss{\'e} class and moreover $(\G_n, \G_n^*)$ is an excellent pair. This can be seen by using a partite method construction as in \cite{NR77,NR83,NR89}, or by introducing new relation for ``not edge".

Let $\T$ be a collection of finite tournaments with $\vert T \vert \geq 3$ for all $T \in \T$, and let $\F(\T)$ denote the class of finite directed graphs $\AA$ with the property that $\BB \not\into \AA$ for all $\BB \in \T$. Then $\F(\T)$ is a Fra{\"i}ss{\'e} class with limit $\FF(\T) := \Flim(\F(\T))$. Let $\F^*(\T)$ denote the class of finite structures of the form $(A, \rightarrow^A, \leq^A)$ with the property that $(A, \rightarrow^A) \in \F(\T)$ and $\leq^A$ is a linear order on $A$. Then $\F^*(\T)$ is a Fra{\"i}ss{\'e} class and moreover $(\F(\T), \F^*(\T))$ is an excellent pair, again see \cite{NR77,NR83,NR89} and \cite{AH78}.

\section{Amenability}
\label{sec:Amen}

Now we are in a position to verify that many of the previously mentioned structures have amenable automorphism groups. After establishing a sufficient density condition for amenability and unique ergodicity, we verify that $\DD_n, \hat{\TT^\omega}, \SS, \GG_n$ and $\FF(\T)$ all have amenable automorphism groups. After establishing Theorem~\ref{dense} verifying that these automorphism groups are amenable will amount to a routine counting of expansions of a structure.

%
\subsection{Density result about amenability.}
\label{sec:dense}

Let $\F$ be a class of finite structures and let $\D \sse \F$. We say that $\D$ is \textbf{cofinal} (or dense) in $\F$ if for every $\AA \in \F$ there is a $\DD \in \D$ such that $\AA \leq \DD$. See section \ref{sec:cons_rand_exp} for notation.

\begin{theo}\label{dense}
Let $(\K,\K^*)$ be an excellent pair of Fra\"iss\'e classes. If for every $\AA \leq \BB$ in $\K$ and every $\expand[\prime]{\AA}, \expand[\prime\prime]{\AA} \in \K^*$ we have
\begin{align}
	\#(\AA^\prime, \BB) = \#(\AA^{\prime\prime}, \BB) \label{ass:samenumberofexp}
\end{align}
then $\Aut(\Flim(\K))$ is amenable. 

Moreover, if there is a cofinal subclass $\D \sse \K$ with the property that for every $\DD \in \D$ and every $\expand[\prime]{\DD}, \expand[\prime\prime]{\DD} \in \K^*$ we have $\expand[\prime]{\DD} \cong \expand[\prime\prime]{\DD}$, then $\Aut(\Flim(\K))$ is uniquely ergodic.
\end{theo}

\proc{Proof.} For $\AA \in \K$ and $\expand[0]{\AA}\in \K^*$ define
\[
	\mu_\AA (\{\AA^0\}) := \frac{1}{\#(\AA)}.
\]

We check that $(\mu_\AA)_{\AA \in \K}$ defines a consistent random expansion. Condition (\textbf{P}) is clear so it remains to check (\textbf{E}).

For $\AA \leq \BB \in \K$ and $\expand[0]{\AA} \in \K^*$ we have
\begin{align*} 
\#(\AA) \cdot  \mu_\AA (\{\AA^0\})  	&= 1 = \sum \{\mu_\BB (\{\BB^*\}) : \expand{\BB}\in\K^*\}\\
																			&= \sum_{\substack{\AA^* \\ \expand{\AA}\in\K^*}} \sum \{ \mu_\BB (\{\BB^*\}) : \expand{\AA} \leq \expand{\BB} \in \K^*\}\\
																			&= \#(\AA) \cdot \sum \{\mu_\BB (\{\BB^*\}) : \expand[0]{\AA} \leq\expand{\BB}\in\K^*\}\\
		\Rightarrow \mu_\AA (\{\AA^0\}) 	&= \sum \{\mu_\BB (\{\BB^*\}) : \expand[0]{\AA} \leq \expand{\BB}\in\K^*\}.
\end{align*}

Here, (\ref{ass:samenumberofexp}) is used for the fourth equality. By Proposition~\ref{threeone} we obtain that $\Aut(\text{Flim}(\K))$ is an amenable group.

In order to show unique ergodicity, it is enough to show uniqueness of a consistent random $\K^*$ expansion by Proposition~\ref{threeone}.

Let $(\mu_\AA)_{\AA \in \K}$ and $(\gamma_\AA)_{\AA\in\K}$ be two consistent random $\K^*$ expansions. From our assumptions and (\textbf{I}) we obtain that $\mu_\DD \equiv \gamma_\DD$ for all $\DD \in \D$. Fix $\AA \in \K$, and find $\DD \in \D$ such that $\AA \leq \DD$. For any fixed $\expand[0]{\AA}\in\K^*$ we have:
\begin{align*}
	\mu_\AA (\{\AA^0\}) 	&\overeasy \sum\{\mu_\DD (\{\DD^*\}) : \expand[0]{\AA} \leq \expand{\DD} \in \K^*\} \\
												&= \sum\{\gamma_\DD (\{\DD^*\}) : \expand[0]{\AA} \leq \expand{\DD} \in \K^*\} \\
												&\overeasy \gamma_\AA (\{\AA^0\})
\end{align*}

Therefore $\mu_\AA \equiv \gamma_\AA$ for all $\AA \in \K$ and uniqueness is verified.
\ep\medbreak

Let us immediately show the usefulness of this result. The class of finite rooted binary trees is not on Cherlin's classification (as it is not a directed graph), but this density result gives amenability very quickly.
 
\subsubsection{Binary trees.}
\label{subsec:binary_trees}

Let $\B$ be the class of finite rooted binary trees. For $\BB \in \B$ we define $T(\BB)$ to be the set of terminal nodes of the tree $\BB$, and we define $\Delta(\BB)$ to be the structure $(T(\BB), C^\BB)$ where $C^\BB$ is a ternary relation on $T(\BB)$. We define $C^\BB$ such that for $x,y,z \in T(\BB)$ we have:
\begin{align*}
	C^\BB(x,y,z) \Leftrightarrow \space & x,y,z \text{ are distinct and the shortest path from } x \text{ to the root}\\ 
																			&\text{is disjoint from the shortest path from } y \text{ to } z.
\end{align*}

In this way we assign to each $\BB \in \B$ a unique $\Delta(\BB)$, but also each structure $\Delta(\BB)$ gives the unique binary tree $\TT$ with the fewest nodes such that $\Delta(\TT) = \Delta(\BB)$. 

Let $\H$ be the class of the structures of the form $\Delta(\BB)$ for $\BB \in \B$. Let $\O\H$ be the class of structures of the form $(A, C^A, \leq^A)$ where $(A,C^A) \in \H$ and $\leq^A$ is a linear ordering of $A$. We say that $\leq^A$ is \textbf{convex} on $(A,\leq^A)$ if
\[
	C^A(x,y,z) \Rightarrow (x <^A y \wedge x<^A z) \vee (y <^A x \wedge z <^A x).
\]

For $\BB \in \B$ and $b \in \BB$ we write $\lev(b) = n$ if the shortest path from $b$ to the root has $n$ edges. We also write $\BB(n) := \{b \in \BB : \lev(b) = n\}$ and $\BB \restrict n$ denotes the subtree given by $\{b \in \BB : \lev(b) \leq n\}$. Define $\BB[b]$ to be the subtree of $\BB$ with root $b$ which contains vertices in $\BB$ whose shortest path to the root of $\BB$ contains $b$. We say that $\BB \in \B$ is an \textbf{$n$-nice tree} if $T(\BB) = \BB(n)$ and $\vert \BB(n) \vert = 2^n$.

We consider the cofinal subclass $\D \sse \H$ which is the collection of structures of the form $(A, C^A)$ for which there exists an $n$-nice tree $\BB \in \B$ such that $\Delta(\BB) = (A, C^A)$. We also have the subclass $\O\D \sse \O\H$, which contains structures in $\D$ augmented by a linear order.

Let $\CH \sse \OH$ that contains the structures with convex linear orderings. Let $\COH$ be the class of structures of the form $(A, C^A, \leq^A, \preceq^A)$ where $(A,C^A,\leq^A) \in \OH$ and $(A,C^A, \preceq^A) \in \CH$. We have that $\H, \OH, \CH$ and $\COH$ are all Fra\"iss\'e classes, see \cite{AN98,BOD12}. Moreover, we have the following.

\begin{cor}\label{cor:OH_amenable} \leavevmode
	\begin{enumerate}
		\item $\Aut(\text{Flim} (\H))$ is uniquely ergodic.
		\item $\Aut(\text{Flim} (\O\H))$ is amenable.
	\end{enumerate}
\end{cor}

\proc{Proof.} \space

(i) We have that $(\H, \CH)$ is an excellent pair, see \cite{MIL79}. We will check the conditions in Theorem~\ref{dense} for the cofinal subclass $\D$.

Fix $\AA \leq \BB$ structures in $\H$ and let $\expand[\prime]{\AA} \in\CH$. Let $\UU \in \B$ be the smallest tree such that $\Delta(\UU) = \BB$. Since this is the smallest tree, each non-terminal node has degree $3$ or $2$. Therefore $\#(\BB) = 2^b$, where $b$ is the number of non-terminal nodes in $\UU$. Similarly, if $\VV$ is the smallest tree such that $\Delta(\VV)=\AA$ then we have $\#(\AA) = 2^a$, where $a$ is the number of non-terminal nodes in $\VV$. Therefore we have that $\#(\AA^\prime,\BB) = 2^{b-a}$ only depends on $\AA$ and $\BB$, not $\AA^\prime$.

(ii) We have that $(\OH, \COH)$ is an excellent pair, see \cite{SOK15}. The conditions in Theorem~\ref{dense}, part 1 follows as in the previous case, but because of rigidity, there is no cofinal class with the isomorphism condition as in part 2. We delay verifying unique ergodicity until Proposition \ref{prop:UE_OH}.

\ep\medbreak

%
\subsection{Generic multipartite digraph $\DD_n$.}

\begin{theo}\label{thm:amen_Dn} For $n \leq \omega$, $\Aut(\DD_n)$ is amenable.
\end{theo}

\proc{Proof.} Since $(\D_n, \D_n^*)$ is an excellent pair, see \cite[Theorem~8.6]{JLNW14} for $n=\omega$ and \cite[Theorem~8.7]{JLNW14} for $n < \omega$, it is enough to show that there is a consistent random $\D_n^*$ expansion of $\D_n$, by Proposition~\ref{threeone}. This will be done using Theorem~\ref{dense}.

Let $\AA = (A, \rightarrow^A)$ and $\BB = (B, \rightarrow^B)$ be structures in $\D_\omega$ such that $\AA \leq \BB$, and let $\leq^A$ be a linear order on $A$ such that $(\AA, \leq^A) \in \D_\omega^*$.

Let $A_1, \ldots, A_a$ be the $\perp^A$-equivalence classes of $\AA$ and let $B_1, \ldots, B_b$ be $\perp^B$-equivalence classes of $\BB$. Without loss of generality we may assume that $\leq^A$ induces a linear ordering on the $\perp^A$-equivalence classes such that $A_1 <^A \ldots <^A A_a$. There are $1 \leq i_1 < i_2 < \ldots < i_a \leq b$ such that $A_j \sse B_{i_j}$ for each $1 \leq j \leq a$.

If $(\BB, \leq^B) \in \D_\omega^*$ is such that $(\AA, \leq^A) \leq (\BB, \leq^B)$, then $\leq^B$ and $\leq^A$ induce the same linear ordering on $\{B_{i_1}, \ldots, B_{i_a}\}$ and they agree on each $A_i$.

Therefore we have:
\begin{align*}
	\#_{\AA,\BB} (\leq^A) 
	&:= \left\vert\{ \leq^B : (\AA, \leq^A) \leq (\BB, \leq^B) \}\right\vert \\
	&= \frac{b!}{a!} \cdot \prod_{k=1}^a \frac{\vert B_{i_k} \vert!}{\vert A_k \vert!} \cdot \prod_{k \notin \{i_1, \ldots, i_a\}} \vert B_k \vert ! \\
	&= \frac{b!}{a!} \cdot \frac{\prod_{k=1}^b \vert B_{i_k} \vert!}{\prod_{k=1}^a \vert A_k \vert!}.
\end{align*}

Moreover, for $\AA_1, \AA_2 \in \binom{\BB}{\AA}$ we have $\#_{\AA_1,\BB} (\leq^{A_1}) = \#_{\AA_2,\BB} (\leq^{A_2})$, so we have that $\Aut(\DD_\omega)$ is amenable.

Now suppose that $\AA$ and $\BB$ are structures in $\D_n$ for a fixed $n < \omega$. We have:
\begin{align*}
	 \#_{\AA,\BB} (\leq^A, \{I_i^A\}_{i=1}^a)
	&:= \left\vert\{ (\leq^B, \{I_i^B\}_{i=1}^b ): (\AA, \leq^A, \{I_i^A\}_{i=1}^a) \leq (\BB, \leq^B, \{I_i^B\}_{i=1}^b) \in \D_n^* \}\right\vert \\
	&= \frac{(n-a)!}{(n-a-b)!} \cdot \prod_{k=1}^a \frac{\vert B_{i_k} \vert!}{\vert A_k \vert!} \cdot \prod_{k \notin \{i_1, \ldots, i_a\}} \vert B_k \vert ! \\
	&= \frac{(n-a)!}{(n-a-b)!} \cdot \frac{\prod_{k=1}^b \vert B_{i_k} \vert!}{\prod_{k=1}^a \vert A_k \vert!}.
\end{align*}

Again, for $\AA_1, \AA_2 \in \binom{\BB}{\AA}$ we have $\#_{\AA_1,\BB} (\leq^{A_1}, \{I_i^{A_1}\}_{i=1}^a) = \#_{\AA_2,\BB} (\leq^{A_2}, \{I_i^{A_2}\}_{i=1}^a)$, so we have that $\Aut(\DD_n)$ is amenable.

\ep\medbreak

%
\subsection{2-cover of the generic tournament $\hat{\TT^\omega}$.}
\label{sec:That_amen}

\begin{theo}\label{thm:amen_T_hat} The group $\Aut(\hat{\TT^\omega})$ is amenable.
\end{theo}

\proc{Proof.} Since $(\hat\T, \hat\T^*)$ is an excellent pair, it is enough to show that there is a consistent random $\hat\T^*$ expansion of $\hat\T$, by Proposition~\ref{threeone}.

For each $\AA \in \hat\T$ we define a measure $\mu_\AA$ by taking:

\[
	\mu_\AA (\{\AA^*\}) := \frac{1}{\# (\AA)}.
\]

We check, using Theorem~\ref{dense}, that this is indeed a random consistent expansion of $\hat\T$. It is enough to show that for $\AA \leq \BB$ in $\hat\T$ and $\AA^*$ with $\expand{\AA} \in \hat\T^*$ the number
\[
	\#_{\hat\T^*} (\AA^*, \BB) := \left\vert \{ \BB^* : \expand{\AA} \leq \expand{\BB} \in \T^* \} \right\vert
\]
depends only on the isomorphism classes of $\AA,\BB$, and notably not on the particular expansion $\AA^*$.

Let $\AA = (A, \rightarrow^A)$ and $\BB = (B, \rightarrow^B)$ be structures in $\hat\T$ such that $\AA \leq \BB$, and let $\leq^A$ be a linear order on $A$ such that $(\AA, \leq^A) \in \hat\T^*$.

Let $A_1, \ldots, A_a$ be the $\perp^A$-equivalence classes of $\AA$ and let $B_1, \ldots, B_b$ be $\perp^B$-equivalence classes of $\BB$. Without loss of generality we may assume that $\leq^A$ induces a linear ordering on the $\perp^A$-equivalence classes such that $A_1 <^A \ldots <^A A_a$. Moreover we may assume that this linear order is induced by $l_i \in A_i$ where $1 \leq i \leq a$. We may also assume that $\vert A_i \vert = 2$ for all $1 \leq i \leq a$.

There are $1 \leq i_1 < i_2 < \ldots < i_a \leq b$ such that $A_j = B_{i_j}$ for each $1 \leq j \leq a$. Define $I := \{i_1, i_2, \ldots, i_a\}$ and $J := \{1,2, \ldots, b\} \setminus I$.

If $\expand{\BB} = (\BB, \leq^B, I_0^B, I_1^B) \in \hat\T^*$ is such that $\expand{\AA} = (\AA, \leq^A, I_0^A, I_1^A) \leq (\BB, \leq^B, I_0^B, I_1^B)$, then $\leq^B$ and $\leq^A$ induce the same linear ordering on $\{B_{i_1}, \ldots, B_{i_a}\}$ and they agree on each $A_i$. Furthermore $I_0^A = I_0^B$ and $I_1^A = I_1^B$ on $A_i$ (for $1 \leq i \leq a$).

Therefore we have:
\[
	\#_{\hat\T^*} (\AA^*,\BB) 
	:= \left\vert \{ \BB^* : \expand{\AA} \leq \expand{\BB} \in \hat\T^* \} \right\vert 
	= \frac{b!}{a!} \cdot 2^{\vert J \vert}.
\]

Clearly this does not depend on the particular expansion $\AA^*$, or the particular embedding of $\expand{\AA}$ into $\expand{\BB}$. Thus by Theorem~\ref{dense}, we have the desired amenability.
\ep\medbreak

Interestingly, this argument does not work for $\Aut(\hat{\QQ})$, despite the superficial similarities between $\hat{\QQ}$ and $\hat{\TT^\omega}$. The rough idea is that $\QQ$ comes with its linear order, and the precompact expansions of $\hat{\QQ}$ must cohere with this linear order. This severely restricts the number of expansions of a finite substructure of $\hat{\QQ}$. With $\hat{\TT^\omega}$ there is no such linear order that must be cohered with. See Theorem~\ref{thm:Qhat} for more details.

This can be extended to a rough heuristic: ``If the expansions $\K^*$ reference a linear order on $\K$, then $\Aut(\Flim(\K))$ is not amenable". This heuristic can be seen in the examples studied in Section~\ref{sec:non-amen}.

%
\subsection{Semi-generic multipartite digraph $\SS$ and variants.}
\label{sec:A_S}

\begin{theo}\label{thm:semi_amen} $\Aut(\SS)$ is amenable.
\end{theo}

\proc{Proof} Since $(\mathcal{S}, \mathcal{S}^*)$ is an excellent pair, see \cite[Lemma~10.7, Lemma~10.8]{JLNW14}, it is enough to show that there is a consistent random $\mathcal{S}^*$ expansion of $\mathcal{S}$, by Proposition~\ref{threeone}.

For each $\AA \in \mathcal{S}$ we define a measure $\mu_\AA$ by taking:

\[
	\mu_\AA (\{\AA^*\}) := \frac{1}{\# (\AA)}.
\]

We check, using Theorem~\ref{dense}, that this is indeed a random consistent expansion of $\mathcal{S}$. It is enough to show that for $\AA \leq \BB$ in $\mathcal{S}$ and $\AA^*$ with $\expand{\AA} \in \mathcal{S}$ the number
\[
	\#_{\mathcal{S}^*} (\AA^*, \BB) := \left\vert \{ \BB^* : \expand{\AA} \leq \expand{\BB} \in \mathcal{S}^* \} \right\vert
\]
depends only on the isomorphism classes of $\AA,\BB$ and $\expand{\AA}$, and notably not on the particular expansion of $\AA$.

Let $A_1, \ldots, A_a$ be the list of $\perp^A$-equivalence classes in $\AA$, and let $\leq^A$ be a linear ordering on $\AA^*$ such that $A_1 \leq^A \ldots \leq^A A_a$, and let $R^A$ be the binary relation on $\AA^*$. Similarly, let $B_1, \ldots, B_b$ be the list of $\perp^B$-equivalence classes in $\BB$, and let $i_1 < i_2 < \ldots < i_a$ be such that $A_j \sse B_{i_j}$ for $1 \leq j \leq a$.

Let $\BB^* = (B, R^B, \leq^B)$ be such that $\expand{\BB} \in \mathcal{S}^*$. If $\expand{\AA} \leq \expand{\BB}$ then $\leq^B$ extends $\leq^A$ and $R^B$ extends $R^A$. So we have:
\begin{align*}
	\#_{\mathcal{S}^*} (\AA^*, \BB) 	&= \frac{b!}{a!} \cdot \prod_{j=1}^a \frac{b_{i_j}!}{a_j !} 
									\cdot \left(\prod_{\substack{1 \leq j \leq b, \\ j \notin \{i_1, \ldots, i_a\}}} b_j !\right) \cdot 2^{\binom{b}{2} - \binom{a}{2}} \\
													&= 	\left( b! \cdot 2^{\binom{b}{2}} \cdot \prod_{j=1}^b b_j! \right) 
															\left( a! \cdot 2^{\binom{a}{2}} \cdot \prod_{j=1}^a a_j! \right)^{-1}.
\end{align*}

Clearly this quotient depends only on the isomorphism classes of $\AA$ and $\BB$.

Suppose that $B_k$ and $B_l$ are $\perp^B$-equivalence classes such that the linear ordering $\leq^B$ implies $B_k < B_l$. Then there are two ways to put $R^B$ between these two classes. Since $R^B$ extends $R^A$ we have to choose $R^B$ among $\binom{b}{2} - \binom{a}{2}$ pairs of column partitions, see Figure~\ref{fig:semi_R} in section 3.6. This is why $2^{\binom{b}{2} - \binom{a}{2}}$ shows up.
\ep\medbreak

\begin{theo}\label{thm:SSexp_amen} $\Aut(\SS_R)$ and $\Aut(\SS_\leq)$ are amenable.
\end{theo}

\proc{Proof} The cases of $\Aut(\SS_R)$ and $\Aut(\SS_\leq)$ are similar to the case of $\Aut(\SS)$, which was just shown.
\ep\medbreak

%
\subsection{Generic omitters $\GG_n$ and $\FF(\T)$.}

Since the expansions of these classes are arbitrary linear orders, $\Aut(\GG_n)$ and $\Aut(\FF(\T))$ are amenable by \cite[Proposition~9.3]{AKL12}. Moreover, the uniform consistent random expansion is indeed a consistent random expansion. By Proposition~\ref{threeone} this ensures amenability of the automorphism groups.

\section{Failures of amenability}
\label{sec:non-amen}

As mentioned at the end of Section~\ref{sec:That_amen}, when the expansions $\K^*$ make reference to a canonical order of $\K$ we expect $\Aut(\Flim(\K))$ to be non-amenable. This heuristic works for the examples we investigate here: $\SS(n), \PP, \PP(3)$ and $\hat{\QQ}$.

%
\subsection{Weak local orders $\SS(n)$.}
\label{sec:S(n)}

The following argument comes from private communication with A. Kechris. The theorem for the case of $n=3$ was first shown in \cite{ZUC13}, although the following argument has a distinct geometric crux.

\begin{theo}\label{thm:Sn_amen} For $n = 2,3$, $\Aut(\SS(n))$ is not amenable.
\end{theo}

\proc{Proof} Fix $n=2$ or $3$. By Proposition~\ref{threeone} it is enough to show that there is no consistent random $\mathcal{S}(n)^*$ expansion of $\mathcal{S}(n)$.

	Suppose for the sake of contradiction that such an expansion $(\mu_\KK)$ exists. Consider the structures $\AA, \BB, \CC \in \mathcal{S}(n)$ where $\AA \leq \BB, \AA \leq \CC$ and $\BB \cap \CC = \AA$. Let:
	\begin{itemize} 
		\item $\AA = (A, \rightarrow^A)$, $A = \{x,y\}$ and   $x \rightarrow^A y$;
		\item $\BB = (B, \rightarrow^B)$, $B = A \cup \{b\}$, $x \rightarrow^B b$ and $b \rightarrow^B y$;
		\item $\CC = (C, \rightarrow^C)$, $C = A \cup \{c\}$, $x \rightarrow^C c$ and $y \rightarrow^C c$.
	\end{itemize}

Consider also expansions (in $\mathcal{S}(n)^*$) of these structures. Let:
	\begin{itemize} 
		\item $\expand{\AA}    	= (\AA, I_0^A,    I_1^A )$,   where $I_0^A = \{x,y\}$;
		\item $\expand{\BB}    	= (\BB, I_0^B,    I_1^B )$,   where $I_0^B = \{x,y,b\}$;
		\item $\expand{\CC}    	= (\CC, I_0^C,    I_1^C )$,   where $I_0^C = \{x,y,c\}$;
		\item $\expand[**]{\CC} = (\CC, I_0^{**}, I_1^{**})$, where $I_0^{**} = \{x,y\}, I_1^{**} =\{c\}$.
	\end{itemize}

\begin{figure}[!ht]
    \centering
    \includegraphics[width=0.8\textwidth]{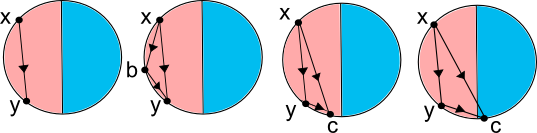}
    \caption{$\expand{\AA}, \expand{\BB}, \expand{\CC}$ and $\expand[**]{\CC}$.}
\end{figure}

This gives us
	\begin{align*}
		\mu_\AA (\{\AA^*\}) &\overeasy \sum\{\mu_\BB (\{\BB^\prime\}) :\expand{\AA} \leq \expand[\prime]{\BB} \in \mathcal{S}(n)^*\} \\
										&= \mu_\BB (\{\BB^*\})
	\end{align*}
since there is only one expansion of $\BB$ in $\mathcal{S}(n)^*$ that extends $\expand{\AA}$, namely $\expand{\BB}$.

Also, we have
	\begin{align*}
		\mu_\AA (\{\AA^*\}) &\overeasy \sum\{\mu_\CC (\{\CC^\prime\}) : \expand{\AA} \leq \expand[\prime]{\CC} \in \mathcal{S}(n)^*\} \\
										&= \mu_\CC (\{\CC^*\}) + \mu_\CC (\{\CC^{**}\})
	\end{align*}
since there are exactly two expansions of $\CC$ in $\mathcal{S}(n)^*$ which extend $\expand{\AA}$, namely $\expand{\CC}$ and $\expand[**]{\CC}$.

Note $\BB \cong \CC$ and this isomorphism extends to all of $\SS(n)$, so $(\textbf{I})$ applies. Also $\expand{\BB} \cong \expand{\CC}$, so we have
\[
		\mu_\CC (\{\CC^*\}) = \mu_\BB (\{\BB^*\}) = \mu_\CC (\{\CC^*\}) + \mu_\CC (\{\CC^{**}\})
\]
So $\mu_\CC (\{\CC^{**}\}) = 0$, which is impossible for a consistent random expansion.
\ep\medbreak

We can also state a related result. Let $\O\mathcal{S}(n)$ be the class of structures of the form $(A, \rightarrow^A, \leq^A)$ where $(A, \rightarrow^A) \in \mathcal{S}(n)$ and $\leq^A$ is a linear order on $A$. Let $\O\mathcal{S}(n)^*$ be the class of the structures of the form $(A, \rightarrow^A, \leq^A, I_0^A, \ldots, I_{n-1}^A)$ where the structures $(A, \rightarrow^A, \leq^A) \in \O\mathcal{S}(n)$ and $(A, \rightarrow^A, I_0^A, \ldots, I_{n-1}^A) \in \mathcal{S}(n)^*$. Then we have that for $n=2,3$, $(\O\mathcal{S}(n), \O\mathcal{S}(n)^*)$ form an excellent pair of Fra\"iss\'e classes which both happen to satisfy the \textbf{SAP}, see \cite{SOK13}. Using the argument in the proof of Theorem~\ref{thm:Sn_amen}, we get the following corollary:

\begin{cor}\label{cor:AutO2} For $n = 2,3$, $\Aut(\Flim(\O\mathcal{S}(n)))$ is not amenable.
\end{cor}

In the case of $n \geq 4$ a similar argument will show that there is a non-trivial flow that is always assigned measure $0$. However, in this case there is no immediate contradiction since the flow need not be minimal; the space might support two disjoint non-minimal flows. A more subtle investigation would be necessary to establish the non-amenability of $\Aut(\SS(n))$, but it eludes the authors at this time. 

%
\subsection{Generic partial order $\PP$.}
\label{sec:P}

The main geometrical idea of the following proof is from \cite[Section 3]{KSOK12}.

\begin{theo} $\Aut(\PP)$ is not amenable.
\end{theo}

\proc{Proof} Since $(\P, \P^*)$ is an excellent pair, it is enough to show that there is no consistent random $\P^*$ expansion of $\P$, by Proposition~\ref{threeone}.

Suppose for the sake of contradiction that such an expansion $(\mu_\KK)$ exists. Consider the structures $\AA, \BB, \CC \in \P$ where $\CC= (C, <^C)$, $C = \{a,b,c\}$, $a <^C b$, with $\AA = \CC \restrict \{a,c\}$ and $\BB = \CC \restrict \{b,c\}$.

Consider also expansions (in $\P^*$) of these structures. Let
	\begin{itemize} 
		\item $\expand{\AA}  = (A, \prec^A)$, where $a \prec^A c$;
		\item $\expand{\BB}  = (B, \prec^B)$, where $b \prec^B c$. 
	\end{itemize}

\begin{figure}[!ht]
    \centering
    \includegraphics[width=0.8\textwidth]{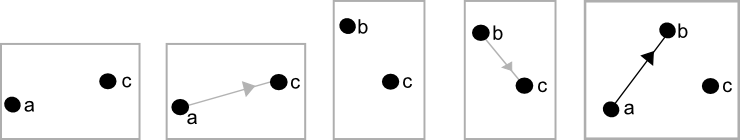}
    \caption{$\AA, \expand{\AA}, \BB, \expand{\BB}$ and $\CC$.}
\end{figure}

Now we consider the possible expansions $\expand{\CC}  = (C, <^C, \prec^C)$ with $\expand{\BB} \leq \expand{\CC}$. We claim that there is only one such expansion possible.

\begin{figure}[!ht]
    \centering
    \includegraphics[width=0.4\textwidth]{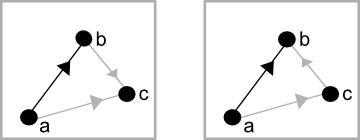}
    \caption{$\expand[0]{\CC}$ and $\expand[1]{\CC}$.}
\end{figure}

We ask ``Is $a \prec^C c$?''. A positive answer yields the expansion $\CC^0$, where $ a \prec^C b \prec^C c$ and $a \prec^C c$, a linear order. A negative answer would yield the cycle $a \prec^C b \prec^C c \prec^C a$, which cannot happen in a linear order. (This is the key geometric observation about $\P^*$.)

Therefore:
	\begin{align*}
		\mu_\BB (\{\BB^*\}) &\overeasy \sum\{\mu_\CC (\{\CC^\prime\}) : \expand{\BB} \leq \expand[\prime]{\CC} \in \P^*\} \\
										&= \mu_\CC (\{\CC^0\}),
	\end{align*}

Now we consider the possible expansions $\expand{\CC}  = (C, <^C, \prec^C)$ with $\expand{\AA} \leq \expand{\CC}$. We claim that there are two such expansions possible.

We ask ``Is $b \prec^C c$?''. A positive answer yields the expansion $\CC^0$. A negative answer yields the expansion $\CC^1$, where $ a \prec^C c \prec^C b$ and $a \prec^C b$, a linear order.

Therefore:
	\begin{align*}
		\mu_\AA (\{\AA^*\}) &\overeasy \sum\{\mu_\CC (\{\CC^\prime\}) : \expand{\AA} \leq \expand[\prime]{\CC} \in \P^*\} \\
										&= \mu_\CC (\{\CC^0\}) + \mu_\CC (\{\CC^1\}).
	\end{align*}

Since $\AA \cong \BB$ and $\expand{\AA} = \expand{\BB}$ we have
\[
	\mu_\CC (\{\CC^0\}) = \mu_\AA (\{\AA^*\}) = \mu_\BB (\{\BB^*\}) = \mu_\CC (\{\CC^0\}) + \mu_\CC (\{\CC^1\})
\]

So $\mu_\CC (\{\CC^1\}) = 0$, which is impossible for a consistent random expansion.
\ep\medbreak

A similar argument will be used to show that $\PP(3)$, the so-called ``twisted generic poset'' has a non-amenable automorphism group. In that case there is an extra layer of notation which somewhat obscures the argument.

%
\subsection{ ``Twisted'' generic partial order $\PP(3)$.}
\label{sec:P(3)}

\begin{theo}\label{thm:AutP3} $\Aut(\PP(3))$ is not amenable.
\end{theo}

\proc{Proof} Since $(\P(3), \P(3)^*)$ is an excellent pair, see \cite[Theorem~9.3]{JLNW14}, it is enough to show that there is no consistent random $\P(3)^*$ expansion of $\P(3)$, by Proposition~\ref{threeone}.

Suppose for the sake of contradiction that such an expansion $(\mu_\KK)$ exists. Consider the structures $\AA, \BB, \CC \in \P(3)$ where $\CC= (C, \leq^C)$, $C = \{a,b,c\}$, $a \rightarrow^C b$, with $\AA = \CC \restrict \{a,c\}$ and $\BB = \CC \restrict \{b,c\}$.

Consider also expansions (in $\P(3)^*$) of these structures. Let
	\begin{itemize} 
		\item $\expand{\AA}  = (A, P_0^A, P_1^A, P_2^A, \preceq^A)$, where $P_0^A = \{a\}$, $P_1^A = \{c\}$ and $a \preceq^A c$;
		\item $\expand{\BB}  = (B, P_0^B, P_1^B, P_2^B, \preceq^B)$, where $P_0^B = \{b\}$, $P_1^B = \{c\}$ and $b \preceq^B c$. 
	\end{itemize}

\begin{figure}[!ht]
    \centering
    \includegraphics[width=0.8\textwidth]{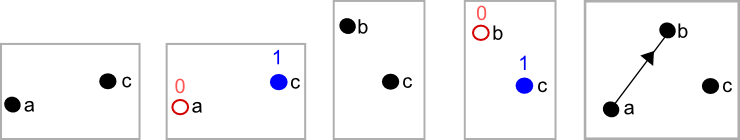}
    \caption{$\AA, \expand{\AA}, \BB, \expand{\BB}$ and $\CC$.}
\end{figure}

Now we consider the possible expansions $\expand{\CC}  = (C, P_0^C, P_1^C, P_2^C, \preceq^C)$ with $\expand{\AA} \leq \expand{\CC}$. There are three options for the label of $b$, namely: $P_0^C(b), P_1^C(b)$ and $P_2^C(b)$. 

\begin{figure}[!ht]
    \centering
    \includegraphics[width=0.8\textwidth]{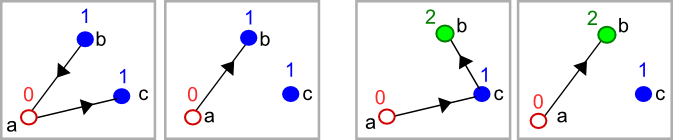}
    \caption{If $P_1^C(b)$, untwisted (1) and twisted (2); If $P_1^C(b)$, untwisted (3) and twisted (4).}
\end{figure}

If $P_1^C (b)$, then we have $b \leq^C a$ and $a \leq^C c$ so $b \leq^C c$. This contradicts the fact that $P_1^C(b)$ and $P_1^C(c)$ guarantee that $b$ and $c$ are $\leq^C$ incomparable.

If $P_2^C (b)$, then we have $a \leq^C c$ and $c \leq^C b$ so $a \leq^C b$. This contradicts the fact that $P_0^C(a), P_2^C(b)$ and $a \rightarrow^C b$ guarantee that $a$ and $b$ are $\leq^C$ incomparable.

Therefore only $P_0^C(b)$ is possible, and so:
	\begin{align*}
		\mu_\AA (\{\AA^*\}) &\overeasy \sum\{\mu_\CC (\{\CC^\prime\}) : \expand{\AA} \leq \expand[\prime]{\CC} \in \P(3)^*\} \\
										&= \mu_\CC (\{\CC^0\}),
	\end{align*}

where $\CC^0$ is given by $P_0^C (b), P_0^C(a), P_1^C (c)$, with $a \rightarrow^C c$ and $a \preceq^C b \preceq^C c$.

On the other hand, there are many expansions $\expand{\CC}$ that respect $\expand{\BB}$, so
	\begin{align*}
		\mu_\BB (\{\BB^*\}) &\overeasy \sum\{\mu_\CC (\{\CC^\prime\}) : \expand{\BB} \leq \expand[\prime]{\CC} \in \P(3)^*\} \\
										&\geq \mu_\CC (\{\CC^0\}) + \mu_\CC (\{\CC^1\}),
	\end{align*}

where $\CC^1$ is given by $P_1^C (a), P_0^C (b), P_1^C(c)$, with $a \rightarrow^C c$ and $b \preceq^C c \preceq^C a$.

\begin{figure}[!ht]
    \centering
    \includegraphics[width=0.8\textwidth]{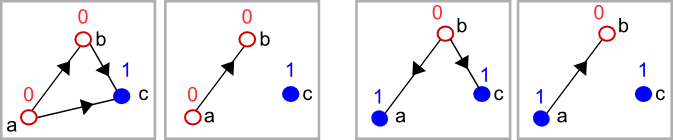}
    \caption{$\expand[0]{\CC}$ untwisted (1) and twisted (2); $\expand[1]{\CC}$ untwisted (3) and twisted (4).}
\end{figure}

Since $\AA \cong \BB$ and $\expand{\AA} \cong \expand{\BB}$ we have
\[
	\mu_\CC (\{\CC^0\}) = \mu_\AA (\{\AA^*\}) = \mu_\BB (\{\BB^*\}) \geq \mu_\CC (\{\CC^0\}) + \mu_\CC (\{\CC^1\}).
\]

So $\mu_\CC (\{\CC^1\}) = 0$, which is impossible for a consistent random expansion.
\ep\medbreak

%
\subsection{2-cover of the linear tournament $\hat{\QQ}$.}
\label{sec:hat_Q}

\begin{theo}\label{thm:Qhat} $\Aut(\hat{\QQ})$ is not amenable.
\end{theo}

\proc{Proof} Let $\T := \Age(\QQ)$.

Suppose for the sake of contradiction that $\Aut(\Flim(\hat{\T}))$ is amenable. By Proposition~\ref{threeone}, there is a consistent random expansion $(\mu_\AA)$. It is enough to work on the cofinal class of structures $\AA$ that have exactly two elements in each $\perp^A$-equivalence class. Moreover, since for all $\expand{\AA}, (\AA\oplus\BB^*) \in \hat{\T}^*$ we have $\expand{\AA} \cong (\AA\oplus\BB^*)$ we must have:
\[
	\mu_\AA (\{\AA^*\}) = \frac{1}{2 \cdot k},
\]
where $k$ is the number of $\perp^A$-equivalence classes, and $\AA = (A, \rightarrow^A)$.

Let $\BB \in \hat{\T}$ be such that each $\perp^B$ equivalence class has two elements with the partition $B = \sqcup_{i=1}^3 B_i$ such that $B_i = \{(i,C), (i,P)\}$. Let the edges on $\{(1,j), (2,j), (3,j)\}$ be given by the natural linear order for $j= C,P$.

Let $\AA$ be the substructure given by the initial segment $B_1 \sqcup B_2$.

Let $\AA^*$ be obtained by the ordering $B_1 <^A B_2$ with $I_1 = \{(1,C),(2,P)\}$ and $I_0 = \{(1,P),(2,C)\}$.

\begin{figure}[!ht]
    \centering
    \includegraphics[scale=0.6]{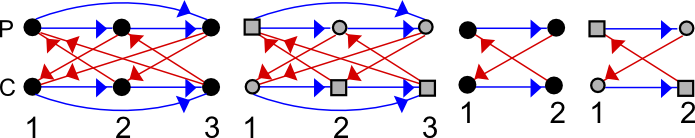}
    \caption{$\BB,\BB^*,\AA,\AA^*$.}
\end{figure}

 Then there is only one $\BB^*$ such that $\expand{\AA} \leq \expand{\BB} \in \hat{\T}^*$, namely $I_1(3,P)$ and $I_0(3,C)$ thus:
\begin{align*}
	\frac{1}{4} = \mu_\AA (\{\AA^*\}) &\overeasy \sum \{\mu_\BB (\{\BB^\prime\}) : \expand{\AA} \leq \expand[\prime]{\BB} \in\hat{\T}^*\} \\
																&= \mu_\BB (\{\BB^*\}) = \frac{1}{6},
\end{align*}
a contradiction.
\ep\medbreak

\section{The Hrushovski property}

We make a brief mention of the Hrushovski property, which is closely related to amenability and unique ergodicity, as examined in \cite{AKL12}. In general, establishing that a Fra\"iss\'e class has the Hrushovski property is a challenging combinatorial problem. Our minor contribution in this area is to establish that the composition $\LEK$ interacts favourably with the Hrushovski property.

\begin{definition} A class $\K$ of finite structures is a \textbf{Hrushovski class} if for any $\KK \in \K$ and any finite sequence of partial isomorphisms $\phi_i : \AA_i \longrightarrow \BB_i$ (for $1 \leq i \leq k$) where $\AA_i, \BB_i \leq \KK$, there is a $\CC \in \K$, such that each $\phi_i$ (for $1 \leq i \leq k$) can be extended to an automorphism $\psi_i : \CC \longrightarrow \CC$.
\end{definition}

Recall the following proposition which appears as Proposition 13.1 in \cite{AKL12} for the special case of order expansions.

\begin{prop}\label{prop:HruAmen} Let $(\K, \K^*)$ be an excellent pair. If $\K$ is a Hrushovski class, then $\Aut(\Flim(\K))$ is amenable.
\end{prop}

We immediately get the following corollary.

\begin{cor} $\mathcal{S}(2)$, $\mathcal{S}(3)$, $\P$, $\P(3)$, $\Age(\SS(2)[\II_n])$, $\Age(\II_n [\SS(2)])$, and $\Age(\hat{\Q})$ are not Hrushovski classes.
\end{cor}

Now we present a proposition which says that the Hrushovski property behaves ``exactly the way you'd want it to'' with respect to the product class $\LEK$.

\begin{prop}\label{prop:HruProd} Let $\K$ and $\L$ be classes of finite relational structures, such that $\L$ satisfies the \textbf{JEP}. Then $\LEK$ is a Hrushovski class if and only if $\L$ and $\K$ are Hrushovski classes.
\end{prop}

\proc{Proof} \space
[$\Rightarrow$] Suppose that $\LEK$ is a Hrushovski class.

First we show that $\L$ is a Hrushovski class. Let $\LL \in \L$, and let $\PP \in \K$ be a one point structure. Fix a finite sequence of partial isomorphisms $\phi_i : \AA_i \longrightarrow \BB_i$ (for $1 \leq i \leq k$) where $\AA_i, \BB_i \leq \LL$. This gives related partial isomorphisms $\phi_i^\prime : (\PP:\AA_i) \longrightarrow (\PP:\BB_i)$. By the Hrushovski property of $\LEK$, there is a $\DD \in \LEK$, and automorphisms $\psi_i^\prime : \DD \longrightarrow \DD$, where $\psi^\prime_i$ extends $\phi_i^\prime$ for $1 \leq i \leq k$. Without loss of generality we may assume that $\DD = (\PP:\CC)$ and consequently, we must have automorphisms $\psi_i : \CC \longrightarrow \CC$ which extends $\psi_i$. So we have verified the Hrushovski property for $\L$.


Now we show that $\K$ is a Hrushovski class. Let $\KK \in \K$, with $\vert \KK \vert = N$ and let $\phi_i : \AA_i \longrightarrow \BB_i$ (for $1 \leq i \leq k$) be a finite sequence of partial isomorphisms, where $\AA_i, \BB_i \leq \KK$. Let $\QQ \in \L$ be a one-point structure. We consider the structure $\DD := (\KK:\underbrace{\QQ, \ldots, \QQ}_{N\text{-times}}) \in \LEK$. Clearly every $\phi_i$ determines a unique $\phi_i^\prime : (\AA_i: \QQ, \ldots, \QQ) \longrightarrow (\BB_i:\QQ, \ldots, \QQ)$, so by the Hrushovski property for $\LEK$ there is a $\CC \in \LEK$ and automorphisms $\psi_i^\prime : \CC \longrightarrow \CC$ where $\psi_i^\prime$ extends $\phi_i^\prime$ for $1 \leq i \leq k$. If $\CC = (\RR:\PP_1, \ldots, \PP_N)$, then every $\phi_i^\prime$ determines an automorphism $\psi_i : \RR \rightarrow \RR$ which extends $\phi_i$.

[$\Leftarrow$] Suppose that $\L$ and $\K$ are Hrushovski classes. Let $\AA = (\KK:\SS_1, \ldots, \SS_N) \in \LEK$ and let $\phi_i : \AA_i \longrightarrow \BB_i$ (for $1 \leq i \leq k$) be a finite sequence of partial isomorphisms in $\AA$. Since $\L$ satisfies the \textbf{JEP} we may assume that $\SS_1 = \ldots = \SS_k =: \SS$.

Each $\phi_i$ is given by a partial isomorphism $\phi_i^\prime : \AA_i^\prime \longrightarrow \BB_i^\prime$ where $\AA_i^\prime, \BB_i^\prime \leq \KK$, and by a sequence of partial isomorphisms $(\phi_{i,s})_{s=1}^{n_i}$ inside $\SS$. Now, by the Hrushovski property for $\K$ there is a $\KK \leq \DD \in \K$ and automorphisms $\psi_i^\prime : \DD \longrightarrow \DD$ which extend the corresponding $\phi_i^\prime$. Moreover, there is an $\SS \leq \TT \in \L$ together with automorphisms $\phi_{i,s}^\prime : \TT \longrightarrow \TT$ which extend the corresponding $\phi_{i,s}^\prime$.

The structure $\EE := (\DD:\TT, \ldots, \TT)$ contains $\AA$ and there are automorphisms $\psi_i : \EE \longrightarrow \EE$ given by $\psi_i^\prime$ and $\psi_{i,s}^\prime$ which extends $\phi_i$. This completes the verification that $\LEK$ is a Hrushovski class.
\ep\medbreak


We are now in a position to give a strengthening of Theorem \ref{A_LEK}.1, with an alternate proof, in the special case that $\L$ and $\K$ are both Hrushovski Classes.

\begin{cor} Let $(\L, \L^*)$ and $(\K, \K^*)$ be excellent pairs of relational structures, where $\L$ and $\K$ are Hrushovski Classes. Then $\Aut(\Flim(\L)), \Aut(\Flim(\K))$ and $\Aut(\Flim(\LEK))$ are amenable.
\end{cor}

\proc{Proof} This follows immediately from Proposition~\ref{prop:HruAmen} and Proposition~\ref{prop:HruProd}.
\ep\medbreak

\section{Unique ergodicity and McDiarmid's inequality}
\label{sec:UE_McD}

This section marks a shift, from establishing amenability to establishing unique ergodicity. Whereas the previous results were exact and finitary, we will now need to make use of asymptotics; we shift from finite combinatorics to analysis.

To establish unique ergodicity we appeal to the probabilistic tools discussed in \cite{AKL12}, suitably generalized to precompact expansions. First we examine the Quantitative Expansion Property, then we will see how this property, together with amenability, gives unique ergodicity. We will compress the probabilistic machinery into the black-box Lemma~\ref{lem:QOP_strategy} which is combinatorial in nature.

%
\subsection{$\QOP$ and $\QOPstar$.}
\label{sec:QOP}

Here we look at two properties that allow us to push amenable automorphism groups up to uniquely ergodic. The following are \textbf{Q}uantitative \textbf{E}xpansion \textbf{P}roperties. The name comes from the property $\mathcal{QOP}$ in \cite{AKL12} which was concerned with expansions that are linear orderings. With suitable adaptions they apply to more general expansions, not just linear orderings.

For fixed structures $\AA$ and $\BB$ with expansions $\expand{\AA}$ and $\expand{\BB}$, and $\E$ a set of embeddings of $\AA$ into $\BB$, define
\[
	\Nexp(\E, \AA^*, \BB^*) := \left\vert\{\phi \in \E : \phi \text{ embeds } \expand{\AA} \text{ into } \expand{\BB}\}\right\vert.
\]
If $\E$ is clear from context, we shall denote this set by $\Nexp(\AA^*,\BB^*)$. Also define 
\[
	\Nemb(\AA,\BB) := \left\vert\{\phi : \phi \text{ embeds } \AA \text{ into } \BB\}\right\vert.
\]
Note that this is $\left\vert \Aut(\AA) \right\vert \cdot \left\vert \binom{\BB}{\AA} \right\vert$, and if $\AA$ is rigid, then this is just $\left\vert \binom{\BB}{\AA} \right\vert$.

\subsubsection{Definitions.}

\begin{definition}[$\QOPstar$] Let $(\K,\K^*)$ be an excellent pair. We say that $\K^*$ satisfies the $\QOPstar$ if there is an isomorphism invariant map $\rho: \K^* \longrightarrow [0,1]$ such that for every $\expand{\AA}\in\K^*$ and every $\epsilon >0$, there is a $\BB \in \K$ and a non empty set of embeddings $\E$, from $\AA$ into $\BB$ with the property that for every $\expand{\BB}\in \K^*$ we have:
\[
\left\vert \frac{\Nexp(\E, \AA^*, \BB^*)}{\vert \E \vert} - \rho\expand{\AA} \right\vert
 < \epsilon.
\]
\end{definition}

Occasionally we will use the notation $a \overset{\epsilon}{\approx} b$ if $\vert a - b \vert < \epsilon$.

\begin{definition}[$\QOP$] Let $(\K,\K^*)$ be an excellent pair. We say that $\K^*$ satisfies the $\QOP$ if there is an isomorphism invariant map $\rho: \K^* \longrightarrow [0,1]$ such that for every $\AA\in\K$ and every $\epsilon >0$, there is a $\BB \in \K$ and a non empty set of embeddings $\E$, from $\AA$ into $\BB$ with the property that for every $\expand{\AA},\expand{\BB}\in \K^*$ we have:
\[
	\left\vert \frac{\Nexp(\E, \AA^*, \BB^*)}{\vert \E \vert} - \rho\expand{\AA} \right\vert
 < \epsilon.
\]
\end{definition}

Note that in general $\QOP$ implies $\QOPstar$ (because $\QOP$ works for an arbitrary expansion, but in $\QOPstar$ you are working with a single expansion). Also, in a Hrushovski class, these are the same (see \cite[Theorem~13.3]{AKL12}), and this is non-trivial.

\subsubsection{General Results and the Main Tool.}

The following theorem is one of the main reasons that we examine $\QOP$. It gives a method for ensuring that an amenable automorphism group is actually uniquely ergodic.

\begin{theo}\label{QOP_UE} Let $(\K,\K^*)$ be an excellent pair, and suppose that $\Aut(\Flim(\K))$ is amenable and $\K^*$ satisfies the $\QOPstar$. Then $\Aut(\Flim(\K))$ is uniquely ergodic.
\end{theo}

\proc{Proof} With minor modifications, this follows from the Fubini-type argument presented in the proof of \cite[Proposition~11.1]{AKL12}.
\ep\medbreak

\subsubsection{Unique Ergodicity of $\Aut(\OH)$.}

Let us illustrate a direct verification of the $\QOPstar$ for a cofinal subclass of $\OH$. Recall the notation from Section \ref{subsec:binary_trees}.

\begin{prop}\label{prop:UE_OH} $\Aut(\Flim(\OH))$ is uniquely ergodic.
\end{prop}

\proc{Proof.} Amenability was proved in Corollary \ref{cor:OH_amenable}.(ii). In order to prove unique ergodicity we will verify $\QOPstar$ for the class $\D$, which is enough by Theorem \ref{dense}.

Let $(A,C^A, \leq^A, \preceq^A) \in \COH$ be given where $(A, C^A, \leq^A) \in \D$, and let $\epsilon > 0$. Let $\BB$ be an $n+1$-nice tree such that $\Delta(\BB) = (A, C^A)$. In particular, we may assume that $A = T(\BB)$. Let $\preceq_1, \ldots, \preceq_l$ be the list of all linear orderings on $A$ such that $(A, C^A, \preceq_i) \in \CH$ for all $1 \leq i \leq l$. Let $\BB^\prime$ be an $n+l$-nice tree such that $\BB \restrict n = \TT$, and let $\BB_1, \ldots, \BB_l$ be the collection of subtrees of $\BB^\prime$ such that for $1 \leq i \leq l$ we have:

\begin{itemize}
	\item $\BB_i \restrict n = \BB$,
	\item $(\forall x \in \BB(n))(\exists!x^\prime \in T(\BB_i))[x^\prime \in \BB^\prime[x]]$,
	\item $T(\BB_i) \sse T(\BB^\prime)$.
\end{itemize}

In this way each $\BB_i$ is a tree of height $n+l$; they are copies of $\BB$ with a branch of height $l$ appended to one of the terminal nodes of $\BB$. Now for $i \neq j$ we have $T(\BB_i) \cap T(\BB_j) = \emptyset$ and $\Delta(\BB_i) \cong \Delta(\BB_j)$.

Moreover, every linear ordering $\preceq^\prime$ such that $(\Delta(\BB_i), \preceq^\prime) \in \CH$ is given by the unique linear ordering with the property that $(\Delta(\TT), \preceq) \in \CH$. More precisely, for $x,y \in T(\BB)$ and $x^\prime \in T(\BB^\prime[x]), y^\prime \in T(\BB^\prime[y])$ we have $x \preceq y \Leftrightarrow x^\prime \preceq^\prime y^\prime$. Therefore, $\preceq^\prime$ is given by one of the $\preceq_i$, and without loss of generality we denote such $\preceq^\prime$ by $\preceq_i$. On each $T(\BB_i)$ we put a linear ordering $\leq_i$ such that
\[
	(T(\BB_i), \leq_i, \preceq_i) \cong (A, C^A, \leq^A, \preceq^A).
\]

Let $\leq$ be a linear ordering on $T(\BB^\prime)$ which extends each $\leq_i$. Note that this is possible since $\BB_i \cap \BB_j = \emptyset$ for $i \neq j$. Let $\phi_i$ be the unique embedding from $(A, C^A, \leq^A)$ into $(\Delta(\BB^\prime), \leq)$ with image $(\Delta(\BB_i), \leq_i)$. Let $\E := \{\phi_i : 1 \leq i \leq l\}$. Now we may take 
\[
	\rho(A, C^A, \leq^A, \preceq^A) := \frac{1}{\#(A, C^A, \leq^A)}.
\]
Now it is easy to see that in this way we can satisfy the condition of the $\QOPstar$, since for a given $\preceq$ with $(\Delta(\BB^\prime), \leq, \preceq)$ there is only one $\phi_i$ such that $\phi_i$ embeds $(A, C^A, \leq^A, \preceq^A)$ into $(\Delta(\BB^\prime), \leq, \preceq)$.
\ep\medbreak

\subsubsection{$\QOP^*, \QOP$ and $\LEK$.}

We now show how the $\QOP$ and the $\QOP^*$ interact with $\LEK$. This will give us an alternate way to check unique ergodicity of $\Aut(\Flim(\LEK))$.

\begin{prop}\label{prop:QOP_LEK} Let $(\K,\K^*)$ and $(\L,\L^*)$ be excellent pairs. If $\K^*$ and $\L^*$ satisfy $\QOPstar$, then $\LEK$ satisfies the $\QOPstar$.
\end{prop}

\proc{Proof.}
Let $\AA = (\SS_1, \ldots, \SS_k : \TT) \in \LEK$. There is an $\RR \in \K$ and an $\E_0$, a collection of embeddings from $\TT$ into $\RR$ which witnesses the $\QOPstar$ for $\K^*$. Also, there are $\LL_i \in \L$ for each $i$ and $\E_i$, a collection of embeddings from $\SS_i$ into $\LL_i$ which witnesses the $\QOPstar$ for $\L^*$.

Consider the structure $\BB = (\LL_1, \ldots, \LL_k : \RR) \in \LEK$ with the collection $\E$ of all embeddings from $\AA$ into $\BB$. Each embedding from $\E$ is given by a member of $\E_0$ and a sequence of embeddings from $\E_1, \ldots, \E_k$.

Let $\mu$ and $\nu$ be maps on $\K$ and $\L$ respectively that verifies the $\QOPstar$. We check that for $\AA^* = (\SS_1^*, \ldots, \SS_k^* : \TT^*)$ with $\expand{\AA} \in \LEKstar$ the following map verifies the $\QOPstar$ for $\LEKstar$:
\[
	\rho\expand{\AA} := \mu\expand{\TT} \cdot \prod_{i=1}^k \nu (\expand{\SS_i})
\]
Notice
\begin{align*}
	\frac{\Nexp(\E, \AA^*, \BB^*)}{\vert \E \vert}
			&= \frac{\Nexp(\E_0, \TT^*, \RR^*) \cdot \prod_{i=1}^k \Nexp(\E_i, \SS_i^*, \LL_i^*)}{\vert \E_0 \vert \cdot \vert \E_1 \vert \cdot \ldots \cdot \vert \E_k \vert}\\
			&\approx (\mu\expand{\TT}\pm \epsilon) \cdot \prod_{i=1}^k (\nu(\expand{\SS_i})\pm \epsilon)\\
			&\overset{(k+1)\epsilon}{\approx} \rho\expand{\AA}
\end{align*}

We trust that the reader can appropriately interpret the use of ``$\approx$'' in the second line. Since $\epsilon$ can be arbitrarily small, this completes the verification of the $\QOPstar$ for $\LEKstar$.
\ep\medbreak

\begin{prop}\label{QOP_LEK} Let $(\K,\K^*)$ and $(\L,\L^*)$ be excellent pairs. If $\K^*$ and $\L^*$ satisfy the $\QOP$, then $\LEK$ satisfies the $\QOP$.
\end{prop}

\begin{cor} Let $(\K,\K^*)$ and $(\L,\L^*)$ be excellent pairs that satisfy the $\QOP$. If $\Aut(\Flim(\L))$ and $\Aut(\Flim(\K))$ are amenable then $\Aut(\Flim(\LEK))$ is uniquely ergodic.
\end{cor}

\proc{Proof.} This follows from Proposition~\ref{ExPair_LEK} (that $(\LEK, \LEKstar)$ is an excellent pair), Proposition~\ref{QOP_LEK}, Theorem~\ref{A_LEK}.i and Theorem~\ref{QOP_UE}.
\ep\medbreak

%
\subsection{Strategy for $\DD_n$ and $\hat{\TT^\omega}$.}
\label{sec:strategyMcD}

Let $\KK$ be one of the directed graphs $\DD_n$ or $\hat{\TT^\omega}$. We will show that $\Aut(\KK)$ is uniquely ergodic using a method developed in \cite[Section 3]{AKL12}. First we present a useful probabilistic inequality, and then we will discuss the general strategy.

\subsubsection{McDiarmid's Inequality.}

The following theorem appears as Lemma 1.2 in \cite{McD89} and is a consequence of Azuma's inequality.

\begin{theo}[McDiarmid's Inequality] Let $\vec{X} = (X_1, \ldots, X_N)$ be a sequence of independent random variables and let $f(X_1, \ldots, X_N)$ be a real-valued function such that there are positive constants $a_i$, with
\[
	\vert f(\vec{X}) - f(\vec{Y}) \vert \leq a_i,
\]

whenever the vectors $\vec{X}$ and $\vec{Y}$ differ only in the $i^{\text{th}}$ coordinate. Then for $\zeta = \EE[f(\vec{X})]$ and all $\epsilon > 0$ we have:
\[
	P[\vert f(\vec{X}) - \zeta \vert \geq \epsilon] \leq 2 \exp\left(\frac{-2\epsilon^2}{\sum_{i=1}^N a_i^2}\right).
\]
\end{theo}

Typically we will use $\vec{X} = (X_1, X_2, \ldots, X_{\binom{n}{2}})$, to talk about structures like the random directed tournament on $n$ vertices, and $\vec{X}$ will correspond to the direction of the $\binom{n}{2}$ edges.

\subsection{General Strategy.}
\label{sec:QOP_strategy}

By Theorem~\ref{QOP_UE}, to show that $\Aut(\KK)$ is uniquely ergodic, it suffices to show that it is amenable and $\K := \Age(\KK)$ satisfies the $\QOPstar$. Showing amenability will usually be direct, and in the case of $\KK = \DD_n$ we have already shown amenability in Theorem \ref{thm:amen_Dn}. Showing that $\K$ satisfies the $\QOPstar$ is a more subtle affair.

For the $\QOPstar$, for a (small) fixed $\HH \in \K$, (with around $k \cdot m$ vertices), we will find a (large, somewhat ``random'') $\GG \in \K$, with $n$ vertices (or sometimes $n$ equivalence classes). To that end, let $\mathbf{G}$ be a uniformly random structure in $\K$ on $n$ fixed vertices (or sometimes $n$ fixed equivalence classes). In general, $\E := \{\phi : \phi \text{ embeds } \AA \text{ into } \BB\}$ will be the set of \textit{all} embeddings from $\HH$ into $\mathbf G$, so $\vert \E \vert =: \Nemb(\HH, \mathbf G)$, and $\rho\expand{\HH} = \frac{1}{\# (\HH)}$, where $\# (\HH)$ is the number of expansions of $\HH$.

We use the notation of $I(n,k,m) := \EE[\Nemb(\HH, \mathbf G)]$, the expected value of the number of embeddings of $\HH$ into $\mathbf G$. Note that $I(n,k,m)$ may also depend on other aspects of $\HH$ and $\mathbf G$, but in practice they won't. In general only $n$ will be allowed to vary, and we will be concerned with large $n$.

We will always establish two separate inequalities using the McDiarmid inequality. The first will be with the function
\[
	f(\GG) := \frac{\Nemb(\HH, \mathbf G)}{I(n,k,m)},
\]
and we will establish that changing $\mathbf G$ by a single edge changes $f(\mathbf G)$ by at most $O(\frac{1}{n^2})$. It is clear that $\EE[f(\mathbf G)] = 1$. McDiarmid's inequality then yields
\begin{align}
	P\left[\left\vert \frac{\Nemb(\HH,\mathbf G)}{I(n,m,k)} - 1 \right\vert \geq D \right]
	&\leq 2 \exp\left( \frac{-2D^2}{\binom{n}{2}\epsilon_1^2 n^{-4}}\right) \nn\\
	&\leq 2 \exp(-\delta_1 D^2 n^2), \label{eqn:f_estimate}
\end{align}
where $D = \frac{\epsilon}{2}$, fixed at the beginning, and $\delta_1$ does not depend on $n$, and the $\epsilon_1$ comes from $O(\frac{1}{n^2})$. The second inequality will be similar, applying McDiarmid's inequality to the function
\[
	f^*(\mathbf G) := \frac{\Nexp(\HH^*,\mathbf{G}^*)}{I(n,k,m)},
\]
and we will establish that changing $\GG$ by a single edge changes $f^*(\mathbf G)$ by at most $O(\frac{1}{n^2})$. It is clear that $\EE[f^*(\mathbf G)] = \rho\expand{\HH}$. Thus McDiarmid's Inequality yields:
\begin{align}
	P\left[\left\vert \frac{\Nexp(\HH^*,\mathbf{G}^*)}{I(n,m,k)} - \rho\expand{\HH} \right\vert \geq D \right] 
	&\leq 2 \exp\left( \frac{-2D^2}{\binom{n}{2}\epsilon_2^2 n^{-4}}\right) \nn\\
	&\leq 2 \exp(-\delta_2 D^2 n^2) \label{eqn:f*_estimate}
\end{align}
where $\rho\expand{\HH}$ and $\delta_2$ do not depend on $n$, and the $\epsilon_2$ comes from $O(\frac{1}{n^2})$.

We then define the probability that $\mathbf G$ is not a suitable candidate:
\[
	p = \frac{2 \cdot \#(\HH)\cdot \#(\mathbf G)}{e^{\delta_3 D^2 n^2}}
\]
which will go to $0$ as $n$ gets large by Stirling's approximation, since for us $\#(\HH) \cdot \#(\mathbf G) = O(n!)$, or $O((n!)^k)$, which corresponds to the number of pairs of expansions on $\HH$ and $\mathbf G$.

So, except with probability $p$, by \ref{eqn:f_estimate} and \ref{eqn:f*_estimate} we have, simultaneously:
\[
	\left\vert \frac{\Nemb(\HH,\mathbf G)}{I(n,m,k)} - 1 \right\vert < D 
	\text{ and }\left\vert \frac{\Nexp(\HH^*,\mathbf{G}^*)}{I(n,m,k)} - \rho\expand{\HH} \right\vert < D
\]
for all $\expand{\HH}, \expand{\mathbf G} \in \K^*$. For large enough $n$, we have $p<1$, so a suitable $\GG$ will exist. The previous inequalities yield the following:
\begin{align*}
	\left\vert \frac{\Nexp(\HH^*,\GG^*)}{I(n,m,k)} - \frac{\Nexp(\HH^*,\GG^*)}{\Nemb(\HH,\GG)} \right\vert 
	&= \frac{\Nexp(\HH^*,\GG^*)}{\Nemb(\HH,\GG)} \cdot \left\vert \frac{\Nemb(\HH,\GG)}{I(n,m,k)} - 1 \right\vert \\
	&\leq 1 \cdot D = D.
\end{align*}
Finally, from the triangle inequality, we have that $\GG$ witnesses the $\QOP$:
\begin{align*}
	 \left\vert \frac{\Nexp(\HH^*,\GG^*)}{\Nemb(\HH,\GG)} - \rho\expand{\HH} \right\vert 
	\leq &\left\vert \frac{\Nexp(\HH^*,\GG^*)}{\Nemb(\HH,\GG)} - \frac{\Nexp(\HH^*,\GG^*)}{I(n,m,k)} \right\vert +
				\left\vert \frac{\Nexp(\HH^*,\GG^*)}{I(n,m,k)} - \rho\expand{\HH} \right\vert \\
\leq &2D \leq \epsilon.
\end{align*}

We summarize this in a lemma.

\begin{lemma}\label{lem:QOP_strategy} Using the notation defined above, suppose that $\Aut(\KK)$ is amenable, that changing $\mathbf G$ by a single edge changes $f$ and $f^*$ by no more than $O(\frac{1}{n^2})$, and that $\#(\mathbf G) \leq O((n!)^k)$. Then $\Aut(\KK)$ is uniquely ergodic.
\end{lemma}

\section{The random method}

We are now in a position to apply Lemma~\ref{lem:QOP_strategy} and establish the unique ergodicity of the automorphism groups of $\DD_\omega, \DD_n, \hat{\TT^\omega}, \SS_\leq$ and $\SS_R$. The structures $\GG_n$ and $\FF(\T)$ are more subtle and require more attention, so they will be addressed in the subsequent section.

%
\subsection{Unique ergodicity of $\Aut(\DD_n)$.}
\label{sec:UEofDn}

We will show unique ergodicity of $\Aut(\DD_n)$ in two steps. First we consider the special case of $n=\omega$, then we adapt the proof for $n < \omega$.

Let $(n)_k$ be the number of injective maps from $\{1, \ldots, k\}$ into $\{1, \ldots, n\}$. Note that in general this is different from $\binom{n}{k}$.

\begin{theo} $\Aut(\DD_\omega)$ is uniquely ergodic.
\end{theo}

\proc{Proof.} Let $\HH = (H, \rightarrow^H) \in \D_\omega$ have $k$ many $\perp^H$-equivalence classes with respective cardinalities $a_1, \ldots, a_k$. Let $G$ be the set with partition $G = \bigsqcup_{i=1}^k G_i$, with $\vert G_i \vert = m$ for $i \leq k$.

We consider a sequence of independent uniformly random variables induced by a pair of elements $G(x,y)$ where $x \in G_i, y \in G_j$ and $i < j$. Each random variable indicates with probability $\frac{1}{2}$ that $x \rightarrow^G y$ and with probability $\frac{1}{2}$ that $y \rightarrow^G x$. In this way, the collection of random variables $G(x,y)$ gives a random directed graph $\mathbf G = (G, \rightarrow^G) \in \D_\omega$.

Notice that $\# (\mathbf G) = (n!)(m!)^n = O(n!)$ and $\# (\HH) = k! \cdot a_1! \cdots a_k!$. We have:
\[
	I(n,m,k,\vec{a}) := \EE[\Nemb(\HH,\mathbf G)] = (n)_k \prod_{i=1}^k (m)_{a_i} \cdot 2^{-\sum_{l<j} a_l a_j},
\]
where $\vec{a} = (a_1, \ldots, a_k)$. For
\[
	f(\mathbf G) := \frac{\Nemb(\HH, \mathbf G)}{I(n,m,k,\vec{a})}
\]
we have $\EE[f(\mathbf G)] = 1$. If we change the direction of only one edge then we change $\Nemb(\HH,\mathbf G)$ by not more than
\[
(k)_2 \cdot (n-2)_{k-2} \cdot \prod_{i=1}^k (m)_{a_i},
\]
and $f$ by not more than:
\[
\frac{(k)_2 \cdot (n-2)_{k-2} \cdot \prod_{i=1}^k (m)_{a_i}}
	{(n)_k \prod_{i=1}^k (m)_{a_i} \cdot 2^{-\sum_{l<j} a_l a_j}} 
	= \frac{1}{n(n-1)} (k)_2 \cdot 2^{\sum_{l<j} a_l a_j} 
	\leq \frac{\epsilon_1}{n^2}
\]
for large enough $n$ and some positive constant $\epsilon_1$.

Let $(\HH, \leq^H) \in \D_\omega^*$ and $(\mathbf G, \leq^G) \in \D_\omega^*$ be such that the $\perp$-equivalence classes are intervals with respect to $\leq^G$. A change in the direction of one edge will change the function:
\[
	f^*(\mathbf G) := \frac{\Nexp(\leq^H,\leq^G)}{I(n,m,k,\vec{a})}
\]
by not more than
\[
	\frac{\binom{k}{2} \cdot (n-2)_{k-2} \cdot \prod_{i=1}^k (m)_{a_i}}
	{(n)_k \cdot \prod_{i=1}^k (m)_{a_i} \cdot 2^{-\sum_{l<j} a_l a_j}} 
	= \frac{1}{n(n-1)}\binom{k}{2}\cdot 2^{\sum_{l<j} a_l a_j} 
	\leq \frac{\epsilon_2}{n^2}
\]
for large enough $n$ and some positive constant $\epsilon_2$. For the McDiarmid inequality we use $\rho\expand{\HH} = \frac{1}{k!a_1! \cdots a_k!}$. Thus we are finished by Lemma~\ref{lem:QOP_strategy}.

\ep\medbreak

\begin{theo} For $n <\omega, \Aut(\DD_n)$ is uniquely ergodic.
\end{theo}

\proc{Proof.} Since it is more natural to let $n$ vary, we will show that $\Aut(\DD_N)$ is uniquely ergodic. For $N$ finite we will consider a similar random directed graph $\mathbf G$ given by the parameter $m$, the cardinality of the parts. In this case, $N$ will be fixed and we will adjust $m$. For the ``small'' structure $\HH$, we use the parameter $k$ for its number of parts.

Notice that $\#(\mathbf G) = (N!)(m!)^N = O((m!)^N)$ and $\#(\HH) = (N)_k \cdot a_1! \cdots a_k!$, where $N$ is fixed, and $m$ can vary. So $\rho\expand{\HH} = \frac{1}{(N)_k\cdot a_1! \cdots a_k!}$.

Let 
\[
	f(\mathbf G) := \frac{\Nemb(\HH,\GG)}{I(N,m,k,\vec{a})}.
\]

A single change in the direction of one edge of parts $i,j$ in $\mathbf G$ will change $\Nemb(\HH,\mathbf G)$ by not more than
\[
	(k)_2 \cdot (N-2)_{k-2} \cdot a_i  \cdot (m-1)_{a_i-1} \cdot a_j  \cdot (m-1)_{a_j-1} \cdot \prod_{l \neq i,j} \binom{m}{a_l},
\]
and $f$ by not more than:
\begin{align*}
	& \frac{(k)_2 \cdot (N-2)_{k-2} \cdot a_i  \cdot (m-1)_{a_i-1} \cdot a_j  \cdot (m-1)_{a_j-1} \cdot \prod_{l \neq i,j} \binom{m}{a_l}}
	{(N)_k \cdot (m)_{a_1} \cdot \ldots \cdot (m)_{a_k} \cdot 2^{-\sum_{i^\prime<j^\prime} a_i^\prime a_j^\prime}} \\
	&= \frac{(k)_2}{N(N-1)}\cdot a_i a_j \frac{1}{m^2}\cdot 2^{\sum_{i^\prime<j^\prime} a_i^\prime a_j^\prime} \\
	&\leq \frac{\epsilon_1}{m^2}
\end{align*}
for a large enough $m$ and some positive constant $\epsilon_1$.

For $\expand{\HH}, \expand{\mathbf G} \in \D_N^*$, for
\[
	f^*(\mathbf G) := \frac{\Nexp(\HH^*,\mathbf G^*)}{I(N,m,k,\vec{a})}
\]
a single change in direction of one edge between parts $i$ and $j$ will change $f^*(\mathbf G)$ by not more than
\begin{align*}
	& \frac{a_i  \cdot (m-1)_{a_i-1} \cdot a_j  \cdot (m-1)_{a_j-1} \cdot \prod_{l \neq i,j} \binom{m}{a_l}}
	{(N)_k \cdot (m)_{a_1} \cdot \ldots \cdot (m)_{a_k} \cdot 2^{\sum_{i^\prime<j^\prime} a_i^\prime a_j^\prime}} \\
	&= \frac{1}{(k)_2}\cdot a_i a_j \frac{1}{m^2}\cdot 2^{\sum_{i^\prime<j^\prime} a_i^\prime a_j^\prime} \\
	&\leq \frac{\epsilon_2}{m^2}
\end{align*}
for large enough $m$ and some positive constant $\epsilon_2$. Thus we are finished by Proposition~\ref{lem:QOP_strategy}.
\ep\medbreak

%
\subsection{The 2-cover $\hat{\TT^\omega}$ has a uniquely ergodic automorphism group.}
\label{sec:hat_T_UE}

Recall the notation from Section \ref{sec:def_blowups}. There we established that $\Aut(\hat{\TT^\omega})$ is amenable.

\begin{theo} The group $\Aut(\hat{\TT^\omega})$ is uniquely ergodic.
\end{theo}

\proc{Proof.} Let $\T := \Age(\TT^\omega)$. Since $(\hat{\T},\hat{\T}^*)$ is an excellent pair, we may use Proposition~\ref{threeone} to establish unique ergodicity. We established that $\Aut(\hat{\TT^\omega})$ is amenable in Theorem \ref{thm:amen_T_hat}.

Let $\HH = (H, \rightarrow^H) \in \hat{\T}$ and let $H = H_1 \sqcup \ldots \sqcup H_k$ be the partition into $\perp^H$-equivalence classes, with $\vert H_i \vert = 2$ for all $i \leq k$.

We consider a sequence of independent random variables $G(i,j)$ with $1 \leq i < j \leq n$. Each random variable indicates with probability $\frac{1}{2}$ that there is an edge between the equivalence classes $G_i$ and $G_j$, so that we obtain a graph in $\T$. Observe that there are only two options since, for a given vertex and equivalence class, there is exactly one in and one out vertex in this class. In this way, by doubling the points and making the canonical edge changes, the collection $G(i,j)$ of random variables gives a directed graph $\mathbf G = (G, \rightarrow^G) \in \hat{\T}$.

Notice that $\# (\mathbf G) = n! \cdot 2^n = O(n!)$ and $\# (\HH) = k! \cdot 2^k$ so $\rho\expand{\HH} = \frac{1}{k! \cdot 2^k}$. In particular we have:
\[
	I(n,k) := \EE[\Nemb(\HH,\mathbf G)] = (n)_k \cdot 2^k \cdot 2^{-\binom{k}{2}} = (n)_k \cdot 2^{-\binom{k}{2} +k}.
\] 
Define
\[
	f(\mathbf G) := \frac{\Nemb(\HH,\mathbf G)}{I(n,k)}
\]
so we have $\EE[f(\mathbf G)] = 1$. Changing a single value of a single $G(i,j)$ changes $\Nemb(\HH,\mathbf G)$ by not more than:
\[
	(k)_2 \cdot (n-2)_{k-2} \cdot 2^k.
\]
and $f(\mathbf G)$ by not more than:
\[
	\frac{(k)_2 \cdot (n-2)_{k-2} \cdot 2^k}{(n)_k \cdot 2^{-\binom{k}{2}+k}} 
	= \frac{1}{n(n-1)}\cdot \frac{(k)_2 \cdot 2}{2^{-\binom{k}{2}-k}} 
	\leq \frac{\epsilon_1}{n^2}
\]
for a large enough $n$ and some positive constant $\epsilon_1$.

Now let $(\HH,\leq^H, I_1^H, I_2^H) \in \hat{\T}^*$ and let $\leq^G, I_1^G, I_2^G$ be given such that $(\mathbf G, \leq^G, I_1^G, I_2^G) \in \hat{\T}^*$. That is, each set $G_i$ comes with a partition given by $I_1^G$ and $I_2^G$, where $\leq^G$ is a linear ordering such that $G_1 \leq^G \ldots \leq^G G_k$. For 
\[
	f^*(\mathbf G) = \frac{\Nexp(\HH^*,\mathbf G^*)}{(n)_k \cdot 2^{-\binom{k}{2} +k}}
\]
a change in a single $G(i,j)$ will change $f^*(\mathbf G)$ by not more than
\[
	\frac{(k)_2 \cdot (n-2)_{k-2}}{(n)_k \cdot 2^{-\binom{k}{2}+k}} = \frac{1}{n(n-1)}\cdot \frac{(k)_2}{2^{-\binom{k}{2}+k}} \leq \frac{\epsilon_2}{n^2}
\]
for a large enough $n$ and a fixed $k$. Thus we are finished by Proposition~\ref{lem:QOP_strategy}.
\ep\medbreak

%
\subsection{Expansions of the semi-generic digraph.}
\label{sec:semi_exp}

In Section \ref{sec:A_S} we established that the automorphism group of $\SS$, the semi-generic digraph, is amenable. In this section we provide some related expansions of $\mathcal{S} := \Age(\SS)$ and check that they satisfy the $\QOP$. This is not enough to get unique ergodicity of $\Aut(\SS)$, but provides a stepping stone to that result, and hones in on the difficulties it presents.

Consider the classes $\mathcal{S}_R := \mathcal{S}^* \vert \{\rightarrow, R\}$ and $\mathcal{S}_\leq := \mathcal{S}^* \vert \{\rightarrow, \leq\}$. It is not hard to see that $\mathcal{S}_R$ and $\mathcal{S}_\leq$ are Fra\"iss\'e classes, and that $(\mathcal{S}_R, \mathcal{S}^*)$ and $(\mathcal{S}_\leq, \mathcal{S}^*)$ are excellent pairs. Let $\SS_R := \Flim(\mathcal{S}_R)$ and $\SS_\leq := \Flim(\mathcal{S}_\leq)$.

\begin{theo}\label{thm:S_exp_UE} $\Aut(\SS_R)$ is uniquely ergodic. 
\end{theo}

\proc{Proof for $\SS_R$.}

Since we already have that $\Aut(\SS_R)$ is amenable, it is enough to prove the uniqueness of a consistent random expansion on a cofinal subclass of $\mathcal{S}_R$.

Let $\HH := (H, \rightarrow^H, R^H) \in \mathcal{S}_R$ and let $H = H_1 \sqcup \ldots \sqcup H_k$ be the partition into $\perp^H$-equivalence classes, with $M =\vert H_i \vert = 2^{k-1}\cdot m$ for all $i \leq k$, for some natural number $m \geq 1$. Moreover, assume that $(H, \rightarrow^H)$ is in the cofinal subclass of $\mathcal{S}$ where all parts in each column have the same size. By part in each column, we mean each element of the partition of a column given by the other columns.

Let $G = \bigsqcup_{i=1}^n G_i$ be a partition with $\vert G_i \vert = M$ for $i \leq n$. We consider a sequence of independent random variables $G(i,j)$ with $1 \leq i < j \leq n$. Each random variable $G(i,j)$ gives a pair of sets $(A,B)$, which is also given by $R$, such that:
\begin{itemize}
	\item $A \sse G_i$,
	\item $B \sse G_j$, and
	\item $\vert A \vert = \vert B \vert = \frac{M}{2}$.
\end{itemize}
The partition given by $(A,B)$ is the same as the partition given by $(G_i \setminus A, G_j \setminus B)$. There are $\frac{1}{2}\cdot \binom{M}{\frac{M}{2}}^2$ such pairs and we assume that $G(i,j)$ has a uniform distribution, i.e. each pair occurs with probability
\[
	p = 2 \cdot \binom{M}{\frac{M}{2}}^{-2}.
\]
Each pair $(A,B)$ describes a distribution of edges between $G_i$ and $G_j$ such that for $x\in G_i, y \in G_j$ we have
\begin{align*}
	x \rightarrow^G y 		&\Leftrightarrow (x\in A, y \in G_j \setminus B) \vee (x \in G_i \setminus A, y \in B)\\
	y \rightarrow^G x 		&\Leftrightarrow (x\in A, y \in B) \vee (x \in G_i \setminus A, y \in G_j \setminus B)\\
	x \in G_i, y \in G_j 	&\Rightarrow (R^G(x,y) \Leftrightarrow y \in B)\\
	x \in G_i, y \in G_j 	&\Rightarrow (R^G(y,x) \Leftrightarrow x \in A)
\end{align*}
In this way we obtain a random structure $\mathbf G = (G, \rightarrow^G, R^G) \in \mathcal{S}_R$. For the McDiarmid Inequality we take $\rho\expand{\HH} = \frac{1}{k! \cdot (M!)^k}$. Notice $\#(\mathbf G) = n!(M!)^n$ and $\#(\HH) = k!(M!)^k$. In particular we have:
\[
	I(n,M,k) := \EE[\Nemb(\HH,\mathbf G)] = (n)_k \cdot (M!)^k \cdot p^{\binom{k}{2}}
\] 
Then for
\[
	f(\mathbf G) := \frac{\Nemb(\HH,\mathbf G)}{I(n,M,k)}
\]
we have $\EE[f(\mathbf G)] = 1$. Changing a single value of a single $G(i,j)$ changes $\Nemb(\HH,\mathbf G)$ by not more than:
\[
	(k)_2 \cdot (n-2)_{k-2} \cdot (M!)^k,
\]
and changes $f(\mathbf G)$ by at most
\[
	\frac{(k)_2 \cdot (n-2)_{k-2} \cdot (M!)^k}{(n)_k \cdot (M!)^k \cdot p^{\binom{k}{2}}} 
	= \frac{1}{n(n-1)}\cdot \frac{(k)_2}{p^{\binom{k}{2}}}
	\leq \frac{\epsilon_1}{n^2}
\]
for a large enough $n$, a positive constant $\epsilon_1$ and a fixed $k$, since $p$ depends only on $M$.

Now let $\expand{\HH},\expand{\mathbf G}\in \mathcal{S}_R$ where $H^* = (H, \leq^H)$ and $\mathbf G^* = (G, \leq^G)$.  Without loss of generality we may assume that $H_1 <^H \ldots <^H H_k$ and $G_1 <^G \ldots <^G G_n$. For 
\[
	f^*(\mathbf G) = \frac{\Nexp(\HH^*,\mathbf G^*)}{I(n,M,k)}
\]
a change in a single $G(i,j)$ will change $f^*(\mathbb G)$ by not more than
\[
	\frac{(k)_2 \cdot (n-2)_{k-2}}{(n)_k \cdot (M!)^k \cdot p^{\binom{k}{2}}}
	= \frac{1}{n(n-1)}\cdot \frac{(k)_2}{ (M!)^k \cdot p^{\binom{k}{2}}}
	\leq \frac{\epsilon_2}{n^2}
\]
for a large enough $n$, a positive constant $\epsilon_2$, and a fixed $k$.  

Thus we are finished by Proposition~\ref{lem:QOP_strategy}.
\ep\medbreak

\subsection{Comments about $\SS$.} The procedure outlined above fails for $\SS$. The major obstacle is that the number of expansions of a structure in $\mathcal{S}$ grows on the order of $O(2^{n^2})$, where $n$ is the number of columns. This invalidates the probabalistic argument presented in section \ref{sec:QOP_strategy}, namely that the probability $p$ of finding a witness $\GG$ does not necessarily limit to $0$.

\section{The hypergraph method}
\label{sec:hypergraph}

In this section we discuss a method for proving the $\QOP$ by using hypergraphs. This method was introduced in \cite{AKL12}, and is different from the one presented in Section~\ref{sec:QOP_strategy}. The idea is to construct a large random object $\GG$ subject to some constraints. We first construct a hypergraph of large girth with many hyperedges. Then each hyperedge is replaced by a random copy of $\HH$. When checking the $\QOP$ for this structure we only examine embeddings that map $\HH$ entirely within a single hyperedge. We shall directly compute the $\QOP$ estimate and will only need a single application of McDiarmid's Inequality.

There is some subtlety in constructing $\GG$ from the hypergraph which is why we include proofs of unique ergodicity of $\Aut(\GG_n)$ and $\Aut(\FF(\T))$, even though similar statement appear in Section 5 of \cite{AKL12}. Our proof of Theorem \ref{thm:UE_hyper} should be compared to the proof of Theorem 5.2 in \cite{AKL12}. 

\subsection{Unique Ergodicity.}

\begin{theo}\label{thm:UE_hyper} Let $n \geq 2$ be a natural number and let $\T$ be a collection of finite tournaments, with $\vert T \vert \geq 3, \forall T \in \T$. Then $\Aut(\GG_n)$ and $\Aut(\FF(\T))$ are uniquely ergodic.
\end{theo}

\proc{Proof.} Since it is more natural to have $n := \vert \GG \vert$ vary, we shall fix $m \geq 2$ and let $n$ vary.

Since $(\F(\T), \F^*(\T))$ and $(\G_m, \G_m^*)$ are both excellent pairs, we may use Proposition~\ref{threeone} to establish amenability and unique ergodicity. Let $\K = \G_m$ or $\F(\T)$, and let $\K^* = \G_m^*$ or $\F(\T)^*$, as appropriate.

Let $\AA$ be a structure in $\K$, and let $\expand{\AA}\in \K^*$. Consider
\[
	\mu_\AA (\{\AA^*\}) := \frac{1}{\vert \AA \vert !}.
\]
By \cite[Proposition~9.3]{AKL12}, $(\mu_\AA)$ is a consistent random expansion since the expansions $\K^*$ of $\K$ are just the linear orders. By Proposition~\ref{threeone} this ensures amenability of $\Aut(\Flim(\K))$.

We check unique ergodicity by verifying the $\QOP$.

Let $\HH \in \K, \expand{\HH} \in \K^*, \vert\HH\vert=k, \rho\expand{\HH} = \frac{1}{k!}$ and $\epsilon >0$. We will find a $\GG \in \K$ and a collection $\E$ of embeddings from $\HH$ into $\GG$ such that for all $\expand{\HH} \in \K^*$ and $\expand{\GG} \in \K^*$ we have:
\[
	\left\vert \frac{\Nexp(\E,\HH^*,\GG^*)}{\vert \E \vert} - \frac{1}{k!} \right\vert < \epsilon.
\]

There is a constant $C$, which depends only on $k$, such that for all $n \geq k$ there is a $k$-uniform hypergraph on $n$ vertices with at least $Cn^{\frac{4}{3}}$ hyperedges and with girth at least $4$, see Lemma 4.1 in \cite{AKL12}. 

Let $n$ be large enough and let $G$ be the underlying set of one such hypergraph, and let $E_1, \ldots, E_s$ be its hyperedges. Since the girth of the hypergraph is at least $4$, for all $i \neq j$ we have in particular that $\vert E_i \cap E_j \vert \leq 1$.

[$\K = \G_m$] Let $x \neq y \in G$. Then:
\begin{itemize}
	\item If there is an $E_i$ such that $\{x,y\} \sse E_i$, there this is exactly one by girth, so then we choose uniformly at random an injective map $e_i: \HH \rightarrow E_i$ and take:
		\[
			x \rightarrow^G y \Leftrightarrow \left(e_i^{-1} (x) \rightarrow^H e_i^{-1} (y) \right).
		\]
	\item Otherwise, fix a directed edge between $x$ and $y$ arbitrarily.
\end{itemize}
In this way we obtain a random directed graph $\mathbf G = (G, \rightarrow^G)$. The construction, and girth at least $4$, ensure that $\II_{m+1}$ can be embedded only in a subgraph induced by $E_i$. However, this is also impossible, since $\HH \in \G_m$, thus $\mathbf G \in \G_m$. In particular, the large girth ensures that $\mathbf G$ does not contain a copy of the three cycle $C_3$.

\begin{figure}[!ht]
    \centering
    \includegraphics[scale=0.6]{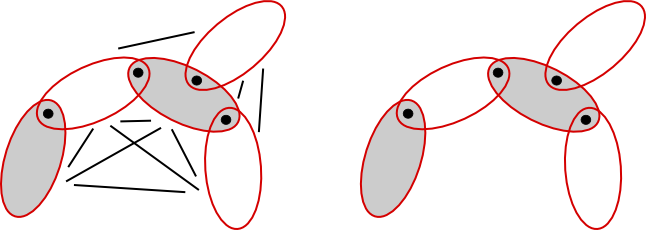}
    \caption{$\GG$ for $\K = \G_m$, and $\GG$ for $\K = \L(\T)$.}
\end{figure}

[$\K = \F(\T)$] Let $x \neq y \in G$. Then:

\begin{itemize}
	\item If there is an $E_i$ such that $\{x,y\} \sse E_i$, there this is exactly one by girth, so then we choose uniformly at random an injective map $e_i: \HH \rightarrow E_i$ and take:
		\[
			x \perp^G y \Leftrightarrow \left(e_i^{-1} (x) \rightarrow^H e_i^{-1} (y) \right).
		\]
	\item Otherwise, $x \perp^G y$.
\end{itemize}
In this way we obtain a random directed graph $\mathbf = (G, \rightarrow^G)$. The construction, and girth at least $4$, ensure that an induced tournament can be embedded only in a subgraph induced by an $E_i$. However, this is also impossible, since $\HH \in \F(\T)$, thus $\mathbf G \in \F(\T)$. In particular, the large girth ensures that $\mathbf G$ does not contain a copy of the three cycle $C_3$.

Now let us check the $\QOP$ estimate. Let $\E$ denote the collection of embeddings of $\HH$ into $\mathbf G$ whose image is completely in one of the $E_i$'s. Note that $\Nemb(\E, \HH,\mathbf G) = s \cdot L$, and $\Nexp (\E, \HH,\mathbf G)$ has a binomial distribution with parameters $(s, \frac{L}{k!})$ where $L = \vert \Aut(\HH) \vert$ and $s \geq Cn^{\frac{4}{3}}$ is the number of hyperedges of $\mathbf G$. Fix structures $\expand{\HH}, \expand{\mathbf G} \in \K^*$.

For
\[
	f(\mathbf G) := \frac{\Nexp (\E, \HH,\mathbf G)}{\Nemb(\E, \HH,\mathbf G)}
\]
we have $\EE[f(\mathbf G)] = \frac{1}{k!}$. Changing a single value of a single $e_i$ changes $f(\mathbf G)$ by not more than:
\[
	\frac{1}{s \cdot L}.
\]
Thus by the McDiarmid inequality we have:
\begin{align*}
	P\left[\left\vert f(\mathbf G) - \frac{1}{k!} \right\vert \geq D \right] 
	&\leq 2 \exp\left( \frac{-2D^2}{s \cdot \left(\frac{1}{s \cdot L}\right)^2}\right) \\
	&\leq 2 \exp\left(-2 \cdot D^2 \cdot L^2 \cdot C\cdot n^{\frac{4}{3}}\right) \\
	&= 2 \exp(-\delta \cdot n^{\frac{4}{3}})
\end{align*}
where $D, L, k$ and $C$ (hence $\delta$) do not depend on $n$. The same estimate holds for all expansions $\expand{\HH}$ and $\expand{\mathbf G}$ in $\K^*$. Therefore, since $\# (\HH) = k!$ and $\# (\mathbf G) = n!$, except on a set of measure
\[
	k! \cdot n! \cdot 2 \exp(-\delta \cdot n^{\frac{4}{3}}),
\]
which is less than $1$ for large $n$, we have
\[
	\left\vert f(\mathbf G) - \frac{1}{k!} \right\vert \leq D.
\]

In particular, choosing $D = \epsilon$ and $n$ large enough, we have our desired digraph $\GG$, which witnesses $\QOP$.
\ep\medbreak

\section{Conclusion and open questions}
\label{sec:conclusion}

The most glaring open question is the following:

\proc{Question 1.}
 Is $\Aut(\SS)$, the automorphism group of the semi-generic multipartite digraph, uniquely ergodic?
\medbreak

Theorem \ref{thm:semi_amen} establishes that it is amenable, and Theorem \ref{thm:S_exp_UE} gives us that the variant $\Aut(\SS_R)$ is uniquely ergodic. It seems as though there are just too many precompact expansions of $\SS$ for the probabilistic methods presented here to work. 

One approach would be to directly analyze the universal minimal flow of $\Aut(\SS)$. One could try a maximal chain construction that was successful for Uspenskij in a related context, see \cite{U00}. See Section 4 of \cite{K12}, and Chapter 6 of \cite{pes06} for good surveys of the results and history relating to the universal minimal flow of the automorphism group of a Fra\"iss\'e structure.

\proc{Question 2.}
 Give a concrete description of the universal minimal flow of $\Aut(\SS)$.
\medbreak

In a separate direction, there are still many open questions relating to Hrushovski classes. The most pressing question here is the following, which was implicitly asked by Hrushovski in \cite{Hru92}.

\proc{Question 3.}
 Is the class of all finite tournaments a Hrushovski class?
\medbreak

This question seems to be quite challenging, given the relatively complicated nature of tournaments. The seemingly easier question about the class $\D_n$ of complete $n$-partite digraphs is also open and still interesting. In general, the known Hrushovski classes seem to all be relational classes with symmetric relations, so any example of a non-symmetric relation Hrushovksi class would be interesting.

\section{Appendix}
\label{sec:appendix}

The appendix contains the proof of the expansion property for $\hat{\T}$.

\subsection{Expansion property for $\hat{\T}$.}

As promised in Section~\ref{sec:def_blowups} we will check that $\hat{\T}^*$ satisfies \textbf{RP} and \textbf{EP}.

\begin{lemma}\label{lem:hatT_Delta} There is a map $\Delta: \hat{\T}^* \longrightarrow \T^*$ which is an injective assignment up to isomorphism, between structures in $\hat{\T}^*$ whose $\perp$-equivalence classes have exactly two elements, and the class $\T^*$ of finite ordered tournaments.
\end{lemma}

\proc{Proof.} Let $(A, \rightarrow^A) \in \hat{\T}$, with $\perp^A$-equivalence classes $A = A_1 \sqcup \ldots \sqcup A_k$ where each class has two elements. Consider a related structure 
\[
	\Delta(\AA) := (\{1,2, \ldots, k\}, \rightarrow^A, \leq^A)
\]
	such that $\leq^A$ is the natural ordering on the set $\{1,2, \ldots, k\}$ and for $1 \leq i < j \leq k$ we have:
\[
	i \rightarrow^A j \Leftrightarrow \left(x \in A_i \wedge y \in A_j \wedge I_1^A(y) \wedge x \rightarrow^A y \right)
\]
Clearly $\Delta(\AA) \in \T^*$. For $\BB = (B, \rightarrow^B, \leq^B) \in \T^*$ we may consider 
\[
	\Delta^{-1}(\BB) := (B \times \{C,P\}, \rightarrow^B, I_0^B, I_1^B, \leq^B) \in \hat{\T}^*
\]
such that for $(x,i), (y,j) \in B \times \{C,P\}$ we have:
\begin{itemize}
	\item $I_1^B ((x,i)) \Leftrightarrow i = C$;
	\item $I_0^B ((x,i)) \Leftrightarrow i = P$;
	\item $(x,i) <^B (y,j) \Leftrightarrow \left( (x=y, i<j) \vee (x <^B y) \right)$;
	\item $(x,i) \rightarrow^B (y,j) \Leftrightarrow \left( (y\rightarrow^B x, i=j) \vee (x \rightarrow^B y, i \neq j) \right)$.
\end{itemize}

\ep\medbreak

Let $\AA^\prime \leq \AA$ be such that $\AA^\prime$ and $\AA$ have the same number of $\perp$-equivalence classes. Then $\left\vert\binom{\AA}{\AA^\prime} \right\vert = 1$. Therefore, in order to verify the \textbf{RP} for $\hat{\T}^*$ it is enough to consider only structures in $\hat{\T}^*$ whose $\perp$-equivalence classes each have exactly two elements.

\begin{theo}\label{thm:hatT_Ramsey} $\hat{\T}^*$ is a Ramsey Class.
\end{theo}

\proc{Proof.} Let $n$ be a natural number and let $\AA,\BB \in \hat{\T}^*$ be such that $\binom{\BB}{\AA} \neq \emptyset$. Without loss of generality we may assume that all $\perp$-equivalence classes in $\AA$ and $\BB$ both have two elements each. Since $\T^*$ is a Ramsey class, see \cite{NR77,NR83,NR89} and \cite{AH78}, there is a (large) $\CC \in \T^*$ such that:
\[
	\CC \longrightarrow \left(\Delta(\BB)\right)^{\Delta(\AA)}_2.
\]
Then we have:
\[
	\Delta^{-1}(\CC) \longrightarrow \left(\BB\right)^{\AA}_2,
\]
and so the verification of the Ramsey Property is complete.
\ep\medbreak

\begin{prop}\label{prop:hatT_EP} $\hat{\T}^*$ satisfies \textbf{EP} with respect to $\hat{\T}$.
\end{prop}

\proc{Proof} We will verify that for each $\AA = (A, \rightarrow^A, \leq^A, I_0^A, I_1^A) \in \hat{\T}^*$ there is an $\HH \in \hat{\T}$ such that for every $\expand{\HH}\in\hat{\T}$ we have $\AA \into \expand{\HH}$.

Since $\hat{\T}^*$ satisfies the \textbf{JEP}, it is enough to obtain \textbf{EP}. Without loss of generality we may assume that each $\perp$-equivalence class in $\AA$ contains exactly two elements.

\begin{figure}[!ht]
    \centering
    \includegraphics[scale=0.6]{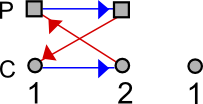}
    \caption{$\XX$ and $\YY$.}
\end{figure}

We make use of the structures $\XX = (X, \rightarrow^X, \leq^X, I_0^X, I_1^X) \in \hat{\T}^*$ and $\YY = (Y, \rightarrow^Y, \leq^Y, I_0^Y, I_1^Y) \in \hat{\T}^*$ such that:
\begin{itemize}
	\item $X = \{1,2\} \times \{P,C\}, Y = \{1\}$,
	\item $I_1^X((1,P)), I_1^X((2,P)), I_1^Y(1)$,
	\item $(1,P) <^X (1,C) <^X (2,P) <^X (2,C)$,
	\item $(1,P) \rightarrow^X (2,P), (1,C) \rightarrow^X (2,C), (2,C) \rightarrow^X (1,P), (2,P) \rightarrow^X (1,C)$.
\end{itemize}

Let $A_1, \ldots, A_k$ be $\perp^A$-equivalence classes which are linearly ordered such that $A_1 <^A A_2 <^A \ldots <^A A_k$. Then there is a $\BB \in \hat{\T}^*$ such that $\AA \leq \BB$ and for every $1 \leq i < k$ there is a $\perp$-equivalence class $B_i^\prime$ in $\BB$ such that:
\begin{enumerate}
	\item $A_i <^B B_i^\prime <^B A_{i+1}$; and
	\item $\BB \restrict (A_i \cup B_i^\prime) \cong \BB \restrict (B_i^\prime \cup A_{i+1}) \cong \XX$.
\end{enumerate}

Then there is a $\BB^\prime \in \hat{\T}^*$ such that $\BB^\prime$ and $\BB$ have the same underlying set and the same relations $I_0, I_1, \rightarrow$ but the linear ordering induced on the $\perp$-equivalence classes in $\BB^\prime$ are opposite to the linear ordering induced on the $\perp$-equivalence classes in $\BB$. Since $\hat{\T}^*$ satisfies the \textbf{JEP} there is a $\CC \in \hat{\T}^*$ such that $\BB \into \CC$ and $\BB^\prime \into \CC$.

Without loss of generality we may assume that each $\perp$-equivalence class in $\CC$ contains exactly two elements. There is a $\CC^\prime \in \hat{\T}^*$ which has the same underlying set as $\CC$, the same linear ordering of equivalence classes, the same $\rightarrow$ relation, but $I_0$ and $I_1$ are inverted. Since $\hat{\T}^*$ satisfies the \textbf{JEP}. there is an $\EE \in \hat{\T}^*$ such that $\CC \into \EE$ and $\CC^\prime \into \EE$. 

\begin{figure}[!ht]
    \centering
    \includegraphics[scale=0.6]{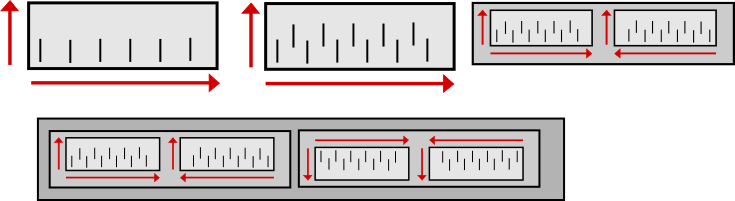}
    \caption{$\AA,\BB,\CC$ and $\EE$.}
\end{figure}

Since $\hat{\T}^*$ is a Ramsey class there are $\FF, \GG \in \hat{\T}^*$ such that:
\[
	\FF \longrightarrow (\EE)^\XX_2 \textnormal{ and } \GG \longrightarrow (\FF)^\YY_2.
\]

Let $\GG = (G, \rightarrow^G, \leq^G, I_0^G, I_1^G)$. We claim that $\HH = (G, \rightarrow^G)$ verifies the \textbf{EP} for $\AA$. Let $\HH^* := (G, \rightarrow^G, \leq^G, I_0^G, I_1^G)$ be such that $\expand{\HH} \in \hat{\T}^*$. Then consider the colouring:
\[
	\xi_Y : \binom{\GG}{\YY} \longrightarrow \{0,1\}
\]
such that
\[
	\xi_Y (\YY^\prime) = 1 \Leftrightarrow I_1^H \restrict Y^\prime = I_1^G \restrict Y^\prime
\]
Consider also the colouring:

\[
	\xi_X : \binom{\GG}{\XX} \longrightarrow \{0,1\}
\]
such that
\begin{align*}
	& \xi_X (\XX^\prime) = 1 \Leftrightarrow \\
	& \leq^H \textnormal{ and } \leq^G \textnormal{ induce the same linear ordering on } \perp\textnormal{-equivalence classes in } \XX^\prime.
\end{align*}

From the construction there are $\FF^\prime \in \binom{\GG}{\FF}$ and $\EE^\prime \in \binom{\FF^\prime}{\EE}$ such that $\xi_Y$ is constantly $c_Y$ on $\binom{\FF^\prime}{\YY}$ and $\xi_X$ is constantly $c_X$ on $\binom{\EE^\prime}{\XX}$.

In particular we have $\xi_Y$ is constant on $\binom{\FF^\prime}{\YY}$. Consider the following options for $(c_X, c_Y)$:

\begin{itemize}
	\item[(1,1)] Here $I_0^H, I_1^H, \leq^H$ agree with $I_0^G, I_1^G, \leq^G$ on $\EE^\prime$ and we have that $\AA \into \EE^\prime$, so $\AA \into \expand{\HH}$.
	\item[(1,0)] Here $I_1^H$ and $I_0^H$ agree with $I_1^G$ and $I_0^G$ on $\BB^\prime$ respectively, but $\leq^H$ and $\leq^G$ induce opposite linear orderings on $\perp$-equivalence classes. Since $\BB^\prime \into \EE^\prime$, this embedding produces $\AA \into \expand{\HH}$.
	\item[(0,1)] Here $I_1^H$ and $I_0^H$ are opposite of $I_1^G$ and $I_0^G$ on $\EE^\prime$ respectively, while $\leq^H$ and $\leq^G$ agree on $\EE^\prime$. Since $\DD^\prime \into \EE^\prime$, this embedding shows that $\AA \into \expand{\HH}$.
	\item[(0,0)] Here $I_1^H, I_0^H, \leq^H$ are opposite of $I_1^G, I_0^G, \leq^G$ on $\EE^\prime$ respectively. Since we have that $\DD^\prime \into \EE^\prime$ and $\BB^\prime \into \EE^\prime$, there is an embedding of $\AA \into \expand{\HH}$.
\end{itemize}

\ep\medbreak


\bibliographystyle{plain}
\bibliography{references}

\end{document}